\documentclass[11pt]{article}
\usepackage{amsmath}
\usepackage{amssymb}

\newcommand{\mm}{\mbox{-}}

\newcommand{\la}{\lambda}

{\catcode `\@=11 \global\let\AddToReset=\@addtoreset}
\AddToReset{equation}{section}

\usepackage{epsfig}
\usepackage{amsmath}
\usepackage{amssymb}
\usepackage{graphicx}
\usepackage{color}

\usepackage[latin1]{inputenc}

\newtheorem{proposition}{Proposition}[section]
\newtheorem{@definition}{\sc Definition}[section]

\newtheorem{@remark}{\sc Remark}[section]
\newenvironment{remark}{\begin{@remark}\rm}{\end{@remark}}
\newtheorem{@example}{\sc Example}[section]

\newcommand{\beqn}{\begin{displaymath}}
\newcommand{\eeqn}{\end{displaymath}}
\newcommand{\beq}{\begin{equation}}  % numbered (single equation)
\newcommand{\eeq}{\end{equation}}

\def\mathsf{\bf}

\def\R{\mathbb{R}}
\def\Z{\mathbb{Z}}

\def\d{\mathrm d}

\def\E{\mathrm E}

\def\text{\mbox}

\def\1{{\bf 1}}

\newcommand{\mbf}[1]{\mbox{\boldmath $#1$}}

\newcommand{\nn}{\nonumber}
\newcommand{\noi}{\noindent}

\newcommand{\mbt}{{\mbf t}}

\newcommand{\mbc}{{\mbf c}}

\newcommand{\mbs}{{\mbf s}}

\newcommand{\mbu}{{\mbf u}}

\newcommand{\mbx}{{\mbf x}}

\newcommand{\mby}{{\mbf y}}

\newcommand{\mba}{{\mbf a}}

\newcommand{\mbq}{{\mbf q}}

\newcommand{\mbgamma}{{\mbf \gamma}}

\newcommand{\smbgamma} {{\small \mbf \gamma}}

\newcommand{\smbt} {{\small \mbf t}}

\newcommand{\smbs} {{\small \mbf s}}

\def\eq2{
\stackrel{\small \rm mod \,2}{=}}

\def\n2{
\stackrel{\small \rm mod \,2}{\neq}}

\def\limd{\renewcommand{\arraystretch}{0.5}
\begin{array}[t]{c}
\stackrel{\rm d}{\longrightarrow} \\
\end{array}\renewcommand{\arraystretch}{1}}

\def\limfdd{\renewcommand{\arraystretch}{0.5}
\begin{array}[t]{c}
\stackrel{\rm fdd}{\longrightarrow} \\
\end{array}\renewcommand{\arraystretch}{1}}

\def\eqfdd{\renewcommand{\arraystretch}{0.5}
\begin{array}[t]{c}
\stackrel{\rm fdd}{=} \\
\end{array}\renewcommand{\arraystretch}{1}}

\def\neqfdd{\renewcommand{\arraystretch}{0.5}
\begin{array}[t]{c}
\stackrel{\rm fdd}{\neq} \\
\end{array}\renewcommand{\arraystretch}{1}}

\newtheorem{thm}{Theorem}[section]
\newtheorem{cor}[thm]{Corollary}

\newtheorem{prop}[thm]{Proposition}

\def\vep{\varepsilon}

\begin{document}

\title{
Anisotropic scaling limits of long-range dependent \\
linear random fields on $\Z^3$
}

%\runtitle{Scaling phase transition for Gaussian fields}

\author{
Donatas Surgailis    \footnote{E-mail: donatas.surgailis@mii.vu.lt}\\
\small \it Vilnius University, Faculty of Mathematics and Informatics, Naugarduko 24, 03225 Vilnius, Lithuania
}
\maketitle

\begin{abstract}

We provide a complete description of anisotropic scaling limits of stationary linear random field
on $\Z^3$ with long-range dependence and moving average coefficients decaying 
as $O(|t_i|^{-q_i})$  in the $i$th direction, 
$i=1,2,3.$ The scaling limits are taken over rectangles in $\Z^3$ whose sides 
increase as $O(\la^{\gamma_i}), i=1,2,3$ when $\la \to \infty$, for any fixed $\gamma_i >0, i=1,2,3 $.
We prove that all these limits are Gaussian RFs whose covariance structure
essentially is determined by the fulfillment or violation
of the balance conditions $\gamma_i q_i = \gamma_j q_j, 1 \le i < j \le 3$. The paper extends
recent results in \cite{ps2015},  \cite{ps2016},   \cite{pils2016},  \cite{pils2017}
on anisotropic scaling  of long-range dependent random fields
from dimension 2 to dimension 3.

\end{abstract}

\smallskip
{\small

\noi {\it Keywords:} scaling  transition; long-range dependence;
linear random field; operator self-similar random field;
fractional Brownian sheet

}

\vskip.7cm

\section{Introduction}

Let $X = \{X(\mbt); \mbt \in \Z^d \}$ be a stationary random field (RF) on $\Z^d, \, d\ge 1, \,
\mbgamma = (\gamma_1, \cdots, \gamma_d) \in \R^d_+ $ be a collection of positive numbers (exponents),
and
\begin{eqnarray} \label{Kp}
K_{\la, \smbgamma} (\mbx)&:=&\prod_{i=1}^d  [1, \lfloor \la^{\gamma_i} x_i\rfloor ] , \qquad  \mbx = (x_1, \cdots, x_d) \in \R^d_+,
%\qquad \mbgamma = (\gamma_1, \dots, \gamma_d) \in \R^d_+
\end{eqnarray}
be a family of $d$-dimensional  `rectangles' indexed by $\la >0$,
whose sides grow at possibly different rate $O(\la^{\gamma_i}), i=1,\cdots, d$ as $\la \to \infty$. Consider
the partial  sums RF:
\begin{eqnarray}\label{Sn}
S^X_{\la, \smbgamma}(\mbx)
&:=&\sum_{{\small \mbt} \in K_{\la, \smbgamma} ({\small \mbx})  }  X(\mbt), \qquad \mbx \in \R^d_+.
\end{eqnarray}
See the end of this section for all
unexplained notation. 
We are interested in the limit distribution of normalized partial sums  \eqref{Sn}:
\begin{equation}\label{partS}
A^{-1}_{\la, \smbgamma} S^X_{\la,\smbgamma}(\mbx)  \ \limfdd \ V^X_{\smbgamma} (\mbx), \quad \mbx \in \R^d_+
\end{equation}
as $\la  \to \infty$,  where $A_{\la, \smbgamma} \to  \infty $ is a normalization. Following \cite{pils2016}, the family
$\{ V^X_{\smbgamma};  \mbgamma \in \R^d_+\} $ of all scaling limits in \eqref{partS} will
be called the {\it scaling diagram of RF $X$}.

The  above problem is classical  for RFs %with long-range dependence (LRD)
except that most previous work dealt with case
$\gamma_1 = \cdots = \gamma_d =1 $ only. See \cite{alb1994}, \cite{dob1979}, \cite{dobmaj1979}, \cite{leo1999}, \cite{leo2011}, 
\cite{sur1982}, \cite{douk2002}, \cite{lav2007}, \cite{lah2016} and the references therein. 
In the latter case,  \eqref{partS}
is naturally referred to as {\it isotropic scaling} while that  %the scaling in  \eqref{partS}
with $\mbgamma \ne (1, \cdots, 1) $  as {\it anisotropic scaling}. For weakly dependent RFs anisotropic scaling
is not very interesting since in such case, summation domains may have very general shape
and the  scaling diagram usually consists of a single point (white noise), or is empty.
See e.g. \cite{bul2012}. Particularly, we note that $d$-dimensional rectangles  in \eqref{Kp} satisfy van Hove's condition
for any $\mbgamma \in \R^d_+ $.

The situation is very different for long-range dependent (LRD) RFs. Although there is no single satisfactory
definition of LRD, usually it refers to stationary RF $X$ with nonsummable covariance function, or unbounded spectral
density, see \cite{dobmaj1979}, \cite{leo1999}, \cite{douk2003p}, \cite{lav2007}, \cite{book2012}.  \cite{ps2016} 
observed that for a large class of LRD RFs $X$ on $\Z^2$,
nontrivial limits
in  \eqref{partS} exist {\it for  any} $\mbgamma = (\gamma_1, \gamma_2) \in \R^2_+$; moreover, there
exists $\gamma^0 >0$ such that
$V^X_{\small \mbgamma} \equiv V^X_\pm $
do not depend on $\mbgamma = (\gamma_1,\gamma_2)$ for
$\gamma_2/\gamma_1 > \gamma^0  $ and $ \gamma_2/\gamma_1 < \gamma^0$, respectively, and the RFs
and $V^X_+, V^X_- $ are different in the sense that
$V^X_+ \neqfdd a V^X_- \, (\forall a >0)$. \cite{ps2016} called
the above phenomenon the {\it scaling transition}. The existence of scaling transition was established in
\cite{ps2015}, \cite{ps2016}, \cite{pils2017}
for a wide class of Gaussian, linear and related nonlinear RFs on $\Z^2 $. It turned out that
for above classes RFs, the scaling limits $V^X_+, V^X_- $ have a very different dependence structure from $V^X_0$, the value
$\gamma^0$ being related to  the intrinsic scale ratio (the ratio of Hurst exponents)
of $X$ along the vertical and horizontal 
axes. Since  $V^X_0, V^X_\pm $ arise in accordance or in violation of the `balance condition' $\gamma_2 = \gamma^0 \gamma_1 $,
\cite{ps2016} termed $ V^X_0$ the {\it well-balanced} and $V^X_\pm $ the {\it unbalanced} scaling limits of $X$, respectively.
A different kind of scaling transition was established for RFs arising by aggregation of network traffic 
and random-coefficient AR(1) time series models in telecommunications and economics, see   \cite{gaig2003},
\cite{miko2002}, \cite{kajt2008}, \cite{pils2014}, \cite{lif2014},  \cite{pils2015}, \cite{lei2018}, also Remark 2.3 
in \cite{ps2016}.  
On the other hand, for some RFs in dimension 2
the scaling diagram may have more than three elements, see  \cite{pils2016}, and
there are classes of LRD RFs which do not exhibit scaling transition (the scaling diagram consists of a single element),
see \cite{ps2015}, \cite{dam2017}.

Since almost all of the above-mentioned work dealt with planar RF models,
a challenging open problem raised in \cite{ps2015}, \cite{pils2017} is anisotropic scaling and identification
of the scaling diagrams of LRD RFs
%extension of scaling transition to RFs
in dimensions 
$d > 2$. The present paper solves  this problem for linear, or moving-average RFs in dimension $d =3 $: % Extension to $d > 3 $ seems feasible, although more cumbersome.
%The  class of RFs for which the above results are established consists of moving-average RFs
\begin{equation}\label{Xlin}
X({\mbf t})\ = \ \sum_{\mbs \in \Z^3} a(\mbt-\mbs) \vep(\mbs), \qquad \mbt = (t_1,t_2,t_3) \in \Z^3,
\end{equation}
where $\{ \vep(\mbs); \mbs  \in \Z^3\}$ is an i.i.d. sequence with zero mean and unit variance, % $\sigma^2>0$,
and $\{a(\mbt), \mbt \in \Z^3\}$ are deterministic
coefficients having the form
\begin{equation} \label{aform}
a(\mbt)\ = \ \frac{g(\mbt)}{ \big(\sum_{j=1}^3 c_j |t_j|_+^{q_j/\nu}% + c_2 |t_2|_+^{q_2/\nu} + c_3 |t_3|_+^{q_3/\nu}
\big)^{\nu}},
\qquad \mbt = (t_1,t_2,t_3) \in \Z^3,
\end{equation}
where $|t|_+ := |t| \vee 1, \, t \in \Z, \,
g(\mbt), \mbt \in \Z^3 $ are bounded with $\lim_{|\mbt| \to \infty} g(\mbt) =: g_\infty \in (0,\infty)$ and
$\nu >0, q_j>0, c_j >0, j=1,2,3$  are parameters satisfying the following inequalities:
\begin{equation}\label{qcond}
1 \ < \ Q := \sum_{i=1}^3 \frac{1}{q_i} \ < \ 2.
\end{equation}
Condition \eqref{qcond} guarantees that
\begin{equation}\label{acond}
\sum_{\mbt \in \Z^3} |a(\mbt)|^2 < \infty  \qquad \text{and} \qquad  \sum_{\mbt \in \Z^3} |a(\mbt)| = \infty.
\end{equation}
In particular,
$X $ in \eqref{Xlin} is a well-defined stationary RF with zero mean, finite variance and covariance
\begin{equation}
\E X({\mbf t}) X({\mbf 0}) = %\sigma^2
\sum_{\mbs \in \Z^3} a(\mbt -  \mbs) a(\mbs), \qquad  \mbt\in \Z^3.
\end{equation}
Notice that $ a(\mbt) = O(|t_i|^{-q_i }) $ as $|t_i| \to \infty $ meaning that when the $q_i$ are different,
the moving-average coefficients decay at different rate
in different directions of $\Z^3$, in which case the RF $X$ exhibits strong anisotropy. On the other hand, when $q_i \equiv q, i=1,2,3$, the RF $X$ is `nearly
isotropic' and conditions \eqref{qcond} reduce to $3/2 < q  < 3$.
%Let us describe the main results of the paper.
The parameter $q_i$ representing `typical scale' of $X$ \eqref{Xlin}
in the $i$th direction, $i=1,2,3$, we may consider
$\gamma^0_{ij} = q_i/q_j,  \, i,j=1,2,3, \, i > j $ as `intrinsic scale ratios' leading to 
three balance conditions
\begin{equation}\label{bal}
\gamma_2/\gamma_1 = \gamma^0_{21}, \qquad \gamma_3/\gamma_1 = \gamma^0_{31}, \qquad  \gamma_3/\gamma_2 = \gamma^0_{32}
\end{equation}
among which only two are independent since any two of \eqref{bal} imply the third one. Depending
on which of the balance conditions in   \eqref{bal} are fulfilled  or violated, we may expect different
scaling limits  $V^X_{\smbgamma}$ of the partial sums
in \eqref{Sn} of the linear RF $X$ \eqref{Xlin}. % which is indeed the case.

%Let us briefly describe
The results of this paper confirm the above intuition.
We prove that for linear RFs in \eqref{Xlin}
the limits $V^X_{\smbgamma} $
in \eqref{partS} exist for any $\mbgamma \in \R^3_+$; moreover depending on $\mbgamma $
these limits can be divided into three groups: the well-balanced limit arising when
$\mbgamma $ satisfies {\it all} balance conditions in \eqref{bal} (group 1);
`partially unbalanced' limits arising when  $\mbgamma $ satisfies {\it only one} of the
balance conditions in \eqref{bal} (group 2), and `completely unbalanced' limits arising when  $\mbgamma $ satisfies {\it none} of the balance conditions in \eqref{bal}
(group 3).
Furthermore, the limit RFs in group 3 agree with FBS (fractional Brownian sheet)  $B_{H_1,H_2,H_3}$ on  $\R^3_+$ with
{\it at least two} among  three indices $H_i \in (0,1], i= 1,2,3 $ equal to 1 or 1/2, while RFs in group 2 are not
FBS but have some `partial FBS property' in one direction.

Let us describe the contents of this work. Sec. 2.1 provides a formal definition
of the  partition of the set $\R^3_+ = \{ \mbgamma \}  $ of scaling exponents into
13 sets $\Gamma_{000}, \cdots, \Gamma_{\mm 1 1 1}$ %indexed by triples,
induced by balance conditions \eqref{bal}. Sec. 2.2 identifies 3 regions (Regions I, II, and III)
in the parameter space
$\{ \mbq = (q_1,q_2,q_3)\in \R^3_+:  1 < Q < 2 \} $ determined by \eqref{qcond} providing
a classification of the linear RF $X$ in \eqref{Xlin} according to the convergence/divergence of
the covariance function on coordinate axes/coordinate planes in $\Z^3 $.
In Sec.~3 we define limit  Gaussian RFs  as stochastic integrals  w.r.t. white noise in $\R^3$  
with kernels taking a different form in Regions I, II, and III, 
and discuss their self-similarity properties. 
We also relate some these limit RFs to FBS % $B_{H_1,H_2,H_3}$
with two Hurst parameters %$H_i, i=1,2,3 $
equal to 1 or $1/2$. Sec.~4 contains the main result (Theorem \ref{main} and Corollary \ref{main1}), 
by identifying {\it all} scaling
limits in \eqref{partS}. Proofs of the main results are given in Sec.~5.

\smallskip

The following comments are in order. We expect that our 
results can be extended for linear RFs in dimension  $d >3 $ with
coefficients $a(\mbt), \mbt \in \Z^d $ having a similar form as in  \eqref{aform} (with \eqref{qcond} replaced
by $1< Q = \sum_{i=1}^d \frac{1}{q_i} < 2 $); however, the description of the scaling limits when
$d > 3 $ seems  cumbersome %and we prefer to stay with the case
and we restrict ourselves  to dimension 
$d=3$ for relative transparency of exposition. Although
the  results
of this paper can be interpreted as a scaling transition occurring at the boundaries of the balance partition, see Fig. 1 below,
we do not attempt to
provide a formal definition of the latter concept for RFs in dimensions $d=3$ or higher. On the other hand,
%since  linear RFs is a particular case of RFs and
%further studies are needed. 
at present
there are many open problems about anisotropic scaling
even for linear RFs in dimension $d=2$. Particularly, we mention the case
%e.g., identification of scaling limits for
(linear) RFs with infinite variance and/or negatively dependent RFs with coefficients as \eqref{aform} but with
$\sum_{\smbt \in \Z^d} |a(\mbt)| < \infty $ (or $Q < 1 $) and
satisfying $\sum_{\smbt \in \Z^d} a(\mbt) = 0 $. See also \cite{lah2016}
on isotropic scaling of negatively dependent linear RFs.

\smallskip

 {\it Notation.} In what follows, $C,  C_1, C_2$ denote  generic positive constants
which may be different at different locations. We write $\limfdd,    \eqfdd, $ and $ \neqfdd $  for the weak convergence, equality and inequality
of finite-dimensional
distributions, respectively.
$\R^d_+ := \{  \mbx =  (x_1, \cdots, x_d) \in \R^d: x_i > 0,
i=1, \cdots, d \},  \bar \R^d_+ := \{ \mbx =  (x_1, \cdots, x_d) \in \R^d: x_i \ge 0,
i=1, \cdots, d \}, \, |\mbx| := \max_{1\le i \le d} |x_i|, 
%\|\mbx\| := (\sum_{i=1}^d |x_i|^2)^{1/2}, \langle \mbx, \mby \rangle := \sum_{i=1}^d x_i y_i, \,
\R_+ := \R^1_+, \bar \R_+ := \bar \R_+, \, \R^2_0 := \R^2 \setminus \{(0,0)\}$.
$\1(A)$ stands for the  indicator function of a set $A$.
%All equalities and inequalities between
%random variables are assumed to hold almost surely.

\section{Preliminaries}

The description of anisotropic limits in \eqref{partS}, or the limiting Gaussian RFs  $V^X_{\smbgamma}$, in the case $d=3$
is considerably more complicated as in the case $d=2$ in
\cite{ps2015}, \cite{pils2017}.  The limit RFs take a different form in different regions depending both on $\mbgamma $ and $\mbq$.
These regions are specified in the following subsections.

\subsection{The balance partition}

For $(\gamma^0_{21}, \gamma^0_{31}, \gamma^0_{32}) \in \R^3_+, \gamma^0_{32} = \gamma^0_{31}/\gamma^0_{21}$
consider the partition
\begin{equation}\label{R3G}
\R^3_+  =  \bigcup_{\imath \in \wp} \Gamma_\imath
\end{equation}
of the set $\R^3_+$ of scaling exponents into 13 sets $\Gamma_\imath, \imath \in \wp $ defined as
\begin{eqnarray*}
\Gamma_{000}
&:=&\{\mbgamma \in \R^3_+:  \gamma_2/\gamma_1 = \gamma^0_{21}, \gamma_3/\gamma_1 = \gamma^0_{31},
\gamma_3/\gamma_2 = \gamma^0_{32}\}, \\
\Gamma_{011}&:=&\{\mbgamma \in \R^3_+:  \gamma_2/\gamma_1 = \gamma^0_{21}, \gamma_3/\gamma_1 > \gamma^0_{31}, \gamma_3/\gamma_2 > \gamma^0_{32} \}, \\
&\dots& \\
\Gamma_{\mm 111}&:=&\{\mbgamma \in \R^3_+:  \gamma_2/\gamma_1 < \gamma^0_{21}, \gamma_3/\gamma_1 > \gamma^0_{31}, \gamma_3/\gamma_2 > \gamma^0_{32}\}.
\end{eqnarray*}
That  is, the index
in $\Gamma_\imath $ is $\imath = k_{21} k_{31} k_{32},  \, k_{ij} \in \{1,0,\mm 1\}$ means that $\gamma_i/\gamma_j > \gamma^0_{ij}$ if $k_{ij} = 1$,
$\gamma_i/\gamma_j = \gamma^0_{ij} $ if $k_{ij} = 0$ and  $\gamma_i/\gamma_j < \gamma^0_{ij}$ if $k_{ij} = \mm 1$, for any  $ 3 \ge i > j \ge 1$.
Thus, the set $\wp = \{\imath\} $ consists of 13 triples:
\begin{equation}
\wp = \{ 000, \ 011, \ 110, \ {10\,\mm 1}, \ {0\,\mm 1\, \mm 1}, \ {\mm 1\mm 10}, \ {\mm 101}, \ {111}, \
 {11\mm 1}, \ {1\mm 1\mm 1}, \ {\mm 1 \mm 1 \mm 1}, \ {\mm 1 \mm 1 1}, \  {\mm 1 1 1} \}.
\end{equation}
The corresponding partition of $\R^2_+ = \{(\gamma_2/\gamma_1, \gamma_3/\gamma_1)\} $ is shown
in Figure 1 below.  There, the line $\Gamma_{000} \subset \R^3_+ $ satisfying all three  balance conditions in \eqref{bal}
(the `well-balanced' set)
reduces to the single point $(\gamma^0_{21}, \gamma^0_{31}) = (q_1/q_2, q_1/q_3)
\in \R^2_+ $, the two-dimensional sets
$\Gamma_{011}, \Gamma_{110}, \Gamma_{10\mm 1}, \Gamma_{0\mm 1\mm 1}, \Gamma_{\mm 1\mm 10}, \Gamma_{\mm 101} $
satisfying only one of the balance conditions in \eqref{bal} (the `partly balanced' sets)
become line segments, and the three-dimensional sets
$\Gamma_{111}, \Gamma_{11\mm 1}, \Gamma_{1\mm 1\mm 1}, \Gamma_{\mm 1 \mm 1 \mm 1}, \Gamma_{\mm 1 \mm 1 1}, \Gamma_{\mm 1 1 1}$
which violate two (or all) balance conditions in \eqref{bal} (the `completely unbalanced' sets) are projected as sets
of dimension 2.

\begin{center}
\begin{figure}[h]
\begin{center}
\includegraphics[width=15 cm,height=20cm]{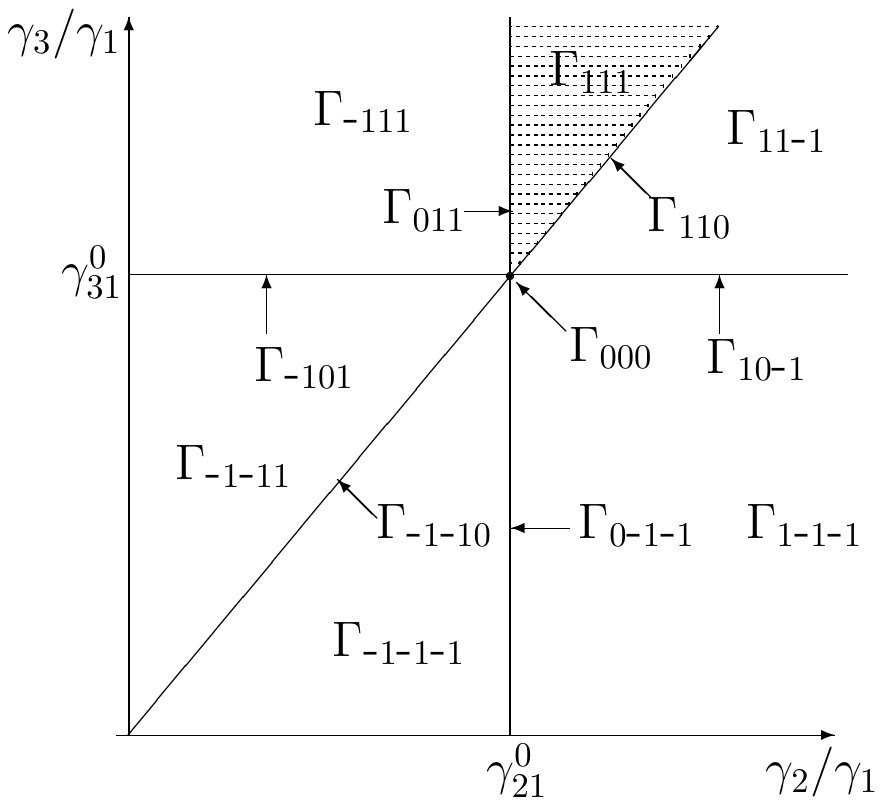}
\vspace{-11.5cm}
\end{center}
\end{figure}
\end{center}

\vskip-2cm

%\center{
\hskip2cm Figure 1. \parbox[t]{12cm} {\small Partition of the quotient space $\R^2_+ = \{(\gamma_2/\gamma_1, \gamma_3/\gamma_1)\} $
induced by balance partition \eqref{R3G}.  The shaded region corresponds to
the subset in \eqref{special}.}

%\vskip1cm

\subsection{Covariance structure of linear LRD RF on $\Z^3$}

As noted above, the scaling  limits  $V^X_{\small \mbgamma}$
depend on parameters $q_j, c_j, j=1,2,3$ in \eqref{aform}. The dependence on the parameters is generally different
in different regions $\Gamma_\imath,  \imath \in \wp$.
Essentially,
it suffices to consider the region
\begin{equation} \label{special}
\{\mbgamma = (\gamma_1,\gamma_2,\gamma_3)\in \R^3_+: \gamma_1 q_1 \le \gamma_2 q_2  \le \gamma_3 q_3\} \ = \ \bigcup_{\imath = 000, 011, 111, 110}
\Gamma_\imath
\end{equation}
(the shaded region in Figure 1) only. Indeed, as shown in Corollary \ref{main1} below, 
for other $\mbgamma$'s,   $V^X_{\small \mbgamma}$ can be defined via a `permutation' of
$q_j, c_j, j=1,2,3$. 
For $\mbgamma $ in \eqref{special},
there are three
parameter regions of $\mbq = (q_1,q_2,q_3)$ defined as follows:
\begin{eqnarray}
&\text{Region I:}\hskip1cm  \frac{1}{2q_1} + \sum_{j=2}^3 \frac{1}{q_j} \ < 1 \ <  \sum_{j=1}^3 \frac{1}{q_j}, \label{R1} \\
&\hskip1cm \text{Region II:}\hskip1cm  \sum_{j=1}^2 \frac{1}{2q_j} + \frac{1}{q_3} < \ 1 \ < \frac{1}{2q_1} + \sum_{j=2}^3 \frac{1}{q_j}, \label{R2}\\
&\text{Region III:}\hskip1cm \sum_{j=1}^3 \frac{1}{2q_j} < 1 < \sum_{j=1}^2 \frac{1}{2q_j} + \frac{1}{q_3}. \label{R3}
\end{eqnarray}
In the `isotropic'
case  $q_1 = q_2 = q_3 =: q$,  \eqref{R1}-\eqref{R3} reduce to $1.5 < q < 2 $ (Region I),
$2 < q < 2.5 $ (Region II) and $2.5 < q < 3 $ (Region III), respectively.

Since the dependence in RF $X$ generally decreases as the $q_j$'s increase, we may say that the dependence in $X$
increases from Region I to Region III. A more precise probabilistic meaning
of the inequalities \eqref{R1}-\eqref{R3} 
in terms of  summability of the covariance function
$r_X (\mbt) := \E X(0) X(\mbt)$
on coordinate axes and coordinate planes in $\Z^3 $
is provided in the following proposition.

\begin{prop}  \label{propcovZ3} Let $X = \{ X(\mbt); \mbt \in \Z^3\}$ be a linear RF in
\eqref{Xlin}-\eqref{aform} % $d=3$ with $q_j >0, j=1,2,3$
satisfying \eqref{qcond}. Then
the covariance  $r_X (\mbt) = \E X(0) X(\mbt), \mbt \in \Z^3$ satisfies the following
properties in respective parameter Regions I-III:
\begin{eqnarray}
\text{Region I:} &&\mbox{$\sum\nolimits_{(t_1,t_2,t_3)\in \Z^3}$}  |r_X(t_1,t_2,t_3)| = \infty, \qquad
\mbox{$\sum\nolimits_{(t_2,t_3)\in \Z^2}$} |r_X(0,t_2,t_3)| < \infty;
\label{R1cov} \\
\text{Region II:} &&\mbox{$\sum\nolimits_{(t_2,t_3)\in \Z^2}$} |r_X(0,t_2,t_3)| = \infty, \qquad
\mbox{$\sum\nolimits_{t_3\in \Z}$} |r_X(0,0,t_3)| < \infty; \label{R2cov}\\
\text{Region III:} &&\mbox{$\sum\nolimits_{t_3\in \Z}$} |r_X(0,0,t_3)| = \infty. \label{R3cov}
\end{eqnarray}
\end{prop}

\begin{remark} The divergence of the series in \eqref{R3cov} can be interpreted as the LRD property
of the sectional process $ \{ X(0,0, t_3); t_3 \in \Z\}$ on the coordinate axis $t_3$ in the parameter Region I.
On the other hand, \eqref{R2cov} say that,  in the parameter Region II
the last process is short-range dependent (SRD)  but the sectional RF  $\{ X(0, t_2,t_3); (t_2,t_3) \in \Z^2 \} $ is LRD. Finally \eqref{R3cov} say
that in the parameter Region III the last sectional RF is SRD  but the RF $X$ on $\Z^3$ is LRD.  
Conditions  \eqref{R1}-\eqref{R3} are not symmetric w.r.t. permutation of $q_j, j=1,2,3 $ and therefore the axes
$t_j, j=1,2,3 $  generally cannot be exchanged in  \eqref{R1cov}-\eqref{R3cov} except for the `isotropic'
case  $q_1 = q_2 = q_3$. 

\end{remark}

\begin{remark} We expect that, under some additional conditions on $g(\mbs)$ in \eqref{acond},
the linear RF $X$ in Proposition
\ref{propcovZ3} has a spectral density % viz., $r_X(\mbt) = \int_{\Pi^3} \e^{\i \langle \smbt, \smbu \rangle} f_X({\mbu}) \d \mbu     $
of the form 
\begin{equation} \label{fspec}
f_X(\mbu) =
\frac{\widetilde g(\mbu )}{\big(\sum_{i=1}^3 \tilde c_i |u_i|^{\alpha_i/\tilde \nu}\big)^{\tilde \nu}}, \qquad \mbu \in \Pi^3 := [-\pi,\pi]^3,
\end{equation}
where $\tilde \nu >0, \tilde c_i >0, \alpha_i>0, i=1,2,3 $ are parameters, $\widetilde g(\mbu ) \ge 0, \mbu \in \Pi^3 $ is a bounded function continuous at the origin with $\widetilde g(\mbf 0) >0$, and the $\alpha_i$'s are related to the $q_i$'s  as
\begin{eqnarray} \label{Hq}
&\alpha_i\ = \ 2q_i\big(\sum_{j=1}^3 \frac{1}{q_j} - 1\big), \qquad q_i \ = \ \alpha_i 
\big(\sum_{j=1}^3 \frac{1}{\alpha_j} - \frac{1}{2}\big). 
\quad i=1,2,3.
\end{eqnarray}
Under \eqref{Hq}, the balance conditions in \eqref{bal} can be rewritten in spectral terms as $\gamma_i/\gamma_j = \alpha_j/\alpha_i, 1\le i <j \le 3$.
See also (\cite{ps2016}, p.2259).  
%Since $\chi_j/\chi_i = q_j/q_i, 1,j=1,2,3 $
Particularly, \eqref{qcond}  is equivalent to $\sum_{i=1}^3 \frac{1}{\alpha_i}  >  1$ and  $\alpha_i >0, i=1,2,3$.
%According to the above conjecture, Regions I-III correspond to 
In terms of `spectral parameters' $\alpha_j, j=1,2,3 $ in \eqref{Hq}, Regions I-III in  \eqref{R1}-\eqref{R3} correspond to
$\alpha_1 <1 $ (Region I), $\frac{1}{\alpha_1} <  1  < \sum_{j=1}^2 \frac{1}{\alpha_j} $ (Region II), and
$\sum_{j=1}^2 \frac{1}{\alpha_j} < 1 < \sum_{j=1}^3 \frac{1}{\alpha_j}$ (Region III).   
The above conjecture agrees with Proposition \ref{propcovZ3}. Indeed, 
the spectral density of the sectional RF  $\{ X(0, t_2,t_3); (t_2,t_3) \in \Z^2 \} $ is
$f_{23}(u_2,u_3) = \int_\Pi f_X (u_1,u_2,u_3) \d u_1 $ which satisfies
$C_1  \bar f_{23}(u_2,u_3)\le f_{23}(u_2,u_3)\le C_2 \bar f_{23}(u_2,u_3),
\bar f_{23}(u_2,u_3) := \int_\Pi  \big(\sum_{i=1}^3  |u_i|^{\alpha_i}\big)^{-1} \d u_1 $, see  \eqref{Cpm}  below.
Clearly, if $\alpha_1 < 1 $ then $\bar f_{23}(u_2,u_3) \le \int_{\Pi} |u_1|^{-\alpha_1} \d u_1 < C$ is bounded
and hence $  f_{23}(u_2,u_3)$ is a bounded function on $\Pi^2$.  The same fact follows from the summability
of the covariance function $r_X (0,t_2,t_3)$ in Region I. Similarly, $1  < \sum_{j=1}^2 \frac{1}{\alpha_j} $
implies that the spectral density $f_{3}(u_3) = \int_{\Pi^2} f_X (u_1,u_2,u_3) \d u_1 \d u_2 $ of the
sectional process $ \{ X(0,0, t_3); t_3 \in \Z\}$ is bounded, which agrees with
the summability
of the covariance function $r_X (0,0,t_3)$ in Region II.

\end{remark}

\noi {\bf Proof of Proposition \ref{propcovZ3}.} We shall use the following elementary inequality.  For any given
$q_j >0, c_j >0, \nu >0, j=1,2,3$ there exist constants  $C_1, C_2 >0$ such that
\begin{equation} \label{Cpm}
C_1 \sum_{j=1}^3 t_j^{q_j} \ \le \
(\sum_{j=1}^3 c_j t_j^{q_j/\nu})^{\nu} \le \ C_2 \sum_{j=1}^3 t_j^{q_j}, \qquad \forall \ \mbt = (t_1,t_2,t_3) \in \R^3_+.
\end{equation}
Indeed, since $c_j t^{q_j/\nu}_j \le (\max_{1\le j \le 3} c_j) (\sum_{j=1}^3 t_j^{q_j})^{1/\nu}, j=1,2,3 $ the second
inequality in \eqref{Cpm} holds with $C_2 =  (3\max_{1\le j \le 3} c_j)^\nu $ and the first inequality
is similar.  Denote
\begin{eqnarray} \label{rhot}
&\rho(\mbt) := \sum_{j=1}^3 |t_j|^{q_j}, %\quad \|\mbt\| := \max_{1 \le j\le 3} |t_j|,
\qquad \mbt \in \R^3,
\end{eqnarray}
then \eqref{Cpm} and \eqref{aform}  imply
\begin{equation} \label{abounds}
C_1 \rho(\mbt)^{-1} \ \le \ |a(\mbt)| \ \le \ C_2 \rho(\mbt)^{-1}, \qquad \mbt \in \Z^3.
\end{equation}
%where $C_i >0, i=1,2 $ are constants independent of $\mbt $.
We claim that \eqref{abounds} imply 
a similar inequality for the covariance $r_X(\mbt) = \sum_{\smbs \in \Z^3} a(\mbs) a(\mbt + \mbs)$, viz.,
\begin{eqnarray} \label{rbounds}
C_1 \rho(\mbt)^{-(2-Q)}(1+ o(1)) \ \le \  |r_X(\mbt)| \ \le \ C_2 \rho(\mbt)^{-(2-Q)}(1+ o(1)),  \qquad  |\mbt| \to \infty,
\end{eqnarray}
where $Q \in (1,2)$ is as in \eqref{qcond}.
%and  $C'_i >0, i=1,2 $ are some constants.
%Q = \sum_{j=1}^3 1/q_j \in (1,2).  $
To check \eqref{rbounds} consider
the convolution
\begin{eqnarray*}
&(\rho^{-1} \star \rho^{-1})(\mbt)
=\int_{\R^3}  \frac{\d s_1 \d s_2 \d s_3}{
(|s_1|^{q_1} + |s_2|^{q_2} +  |s_3|^{q_3})
(|t_1+s_1|^{q_1} + |t_2+s_2|^{q_2} +  |t_3+ s_3|^{q_3})}, \quad \mbt \in \R^3.
\end{eqnarray*}
By change of variables $s_j \to  \rho^{1/q_j} s_j, j=1,2,3 $ we obtain
\begin{eqnarray}\label{rhoQ}
(\rho^{-1} \star \rho^{-1})(\mbt)
&=&L(\mbt)\rho(\mbt)^{-(2-Q)},
\end{eqnarray}
where
\begin{eqnarray} \label{Lang}
L(\mbt):=
\int_{\R^3}  \frac{\d s_1 \d s_2 \d s_3}{
\sum_{j=1}^3 |s_j|^{q_j} \sum_{k=1}^3
|s_k+ \frac{t_k}{\rho(\smbt)^{1/q_k}}|^{q_k}}
= \int_{\R^3}  \frac{\d \smbs}{
\rho(\smbs) \rho(\smbs + \bar \smbt)},   \quad \bar \mbt := \big(\frac{t_1}{\rho(\smbt)^{1/q_1}}, \frac{t_2}{\rho(\smbt)^{1/q_2}},
\frac{t_3}{\rho(\smbt)^{1/q_3}}\big).
 \end{eqnarray}
Let us show that
\begin{eqnarray}\label{Lbounded}
0 \ < \ C_1 \ \le \ L(\mbt) \ \le C_2 \  < \ \infty, \qquad \mbt \in \R^3.
\end{eqnarray}
Let $%\|\mbt \| := \sum_{j=1}^3 |t_j|, \
B_\delta (\mbt) := \{s \in \R^3: \|\mbt - \mbs\| \le \delta \}, \
B_\delta^c(\mbt) := \R^3 \setminus  B_\delta (\mbt)$. Let us prove first that for any $h>0, \delta >0$
\begin{eqnarray}
\int_{B_\delta (\bf 0)}\rho(\mbt)^{-h} \d \mbt &<&\infty   \quad  \Longleftrightarrow \quad Q > h, \label{Qh1} \\
\int_{B^c_\delta (\bf 0)} \rho(\mbt)^{-h} \d \mbt &<&\infty   \quad  \Longleftrightarrow \quad Q < h.  \label{Qh2}
\end{eqnarray}
By \eqref{Cpm},
it suffices to prove \eqref{Qh1}-\eqref{Qh2} for $h = \delta =1 $. We shall often use the elementary inequalities:
\begin{eqnarray}\label{elem}
\hskip-1cm \int_0^1 \frac{\d u}
{(v+ u^q)^h}&\le&C \begin{cases}1, &0<q< 1/h, 0< v < 1, \\
|\log v|, &q=1/h, 0< v < 1, \\
v^{(1/q) -h}, &q >1/h, 0 < v < 1, \\
v^{-h}, &q >0, v \ge 1,
\end{cases} \qquad
\int_0^\infty \frac{\d u}
{(v+ u^q)^h} \le C v^{(1/q) -h}, \ q >1/h.
\end{eqnarray}
Let us prove the converse implication in \eqref{Qh1}, or
%\noi {\it Proof of \eqref{Qh1}.} It suffices
%to prove
$I := \int_{[0,1]^3} \rho(\mbt)^{-1} \d \mbt <\infty $ if $Q > 1 $.
Using  \eqref{elem}, we get
$I < \infty $ when $q_3 < 1 $ and $I \le C \int_0^1 \d t_1 \int_0^1 (t_1^{q_1} + t_2^{q_2})^{(1/q_3)-1} \d t_2 $ when $q_3 > 1 $.
Using \eqref{elem} again we get $I < \infty $ if $q_2 (1 - (1/q_3)) < 1 $ and
$I \le C \int_0^1 t_1^{(q_1/q_2) + (q_1/q_3) - q_1} \d t_1 < \infty $ if $q_2 (1 - (1/q_3)) > 1 $ where the last integral converges
since $Q> 1 $. The case when $q_3 =1 $ and/or $q_2 (1 - (1/q_3)) = 1 $ follow similarly. The remaining
implications in \eqref{Qh1}-\eqref{Qh2} follow in a similar fashion
and we omit the details.
%{\bf [Indeed let use prove the converse implication in \eqref{Qh2} or
%$J:= \int_{\R^3\setminus [0,1]^3} \rho(\mbt)^{-1} \d \mbt <\infty $ if $Q < 1 $. We have
%$J= J_1 + J_2 + J_3 $ where $J_1 = \int_{[0,1]^2 \times [1, \infty)} \cdots,
%J_2 = \int_{[0,1] \times [1, \infty)^2} \cdots, J_3 = \int_{[1, \infty)^3} \cdots $.
%Since $Q < 1 $ implies $q_3 > 1, 1/q_2 + 1/q_3 < 1 $ we get trivially
%$J_1 < \infty, J_2 < \infty $. Finally, using twice the second ineq. in \eqref{elem} we get
%$J_3 \le C \int_{[1,\infty)^2} \d t_1 \d t_2  (t_1^{q_1} + t_2^{q_2})^{1/q_3 - 1}
%\le C \int_{[1,\infty)^2} \d t_1 (t_1^{q_1})^{(1/q_2) - (1- (1/q_3))} < \infty  $ since $Q < 1 $. ]}

\smallskip

Next, we prove \eqref{Lbounded}.
Note  $\rho (\bar \mbt) = 1 $ and therefore $|\bar \mbt | > \delta
\ \forall \mbt \in \R^3$ for some $\delta >0$.
Split $L(\mbt) = L_1(\mbt) + L_1(\mbt)+ L_{12}(\mbt)$, where
$L_1(\mbt) := \int_{B_\delta (0)}
\rho(\mbs)^{-1} \rho(\mbs + \bar \mbt)^{-1} \d \mbs,
L_2(\mbt) := \int_{B_\delta (-\bar \smbt)}
\rho(\mbs)^{-1} \rho(\mbs + \bar \mbt)^{-1} \d \mbs,
L_{12}(\mbt) := \int_{\R^3 \setminus (B_\delta (0) \cup B_\delta (-\bar \smbt))}
\rho(\mbs)^{-1} \rho(\mbs + \bar \mbt)^{-1} \d \mbs $. Since
$\rho(\mbs + \bar \mbt)^{-1} $ is bounded on $B_\delta (0)$ for $\delta >0$ small enough
it follows that $L_1(\mbt) \le C\int_{B_\delta (0)}
\rho(\mbs)^{-1} \d \mbs \le C $ in view of  \eqref{Qh1} and $Q > 1 $. Similarly,
$L_2(\mbt) \le  \int_{B_\delta (-\bar \smbt)}  \rho(\mbs + \bar \mbt)^{-1} \d \mbs \le  C $.
Finally, $|L_{12}(\mbt)|^2 \le
\big(\int_{B^c_\delta(0)} \rho(\mbs)^{-2} \d \mbs\big)^{1/2}
\big(\int_{B^c_\delta(-\bar \smbt)} \rho(\mbs + \bar \mbt)^{-2} \d \mbs\big)^{1/2} \le C $
according to \eqref{Qh1} and $Q < 2 $.  This proves the upper bound in \eqref{Lbounded}.
The lower bound in \eqref{Lbounded}  follows from the uniform boundedness from below
of the integrand of \eqref{Lang} in a vicinity of the origin, viz.,
$\inf_{\bar \smbt \in \R^3} \inf_{\smbs \in B_\delta (0)} \rho(\mbs)^{-1} \rho(\mbs + \bar \mbt)^{-1}  >C >0$ for
any $\delta >0$ small enough. % implying $L_1(\mbt) > C >0 \ \forall \mbt \in \R^3$.

\smallskip

Let us prove \eqref{rbounds}. We use the following inequality: for all
$K > 0$ large enough
\begin{equation}\label{rhobdd}
\rho({\mbs}_1)/2  <  \rho({\mbs}_2) < 2\rho({\mbs}_1),  \quad  \forall \  |{\mbs}_i| > K, \ i=1,2, \   |{\mbs}_1 - {\mbs}_2| \le 1,
\end{equation}
which follows by Taylor expansion of  $\rho({\mbt})$ in \eqref{rhot}.
For a large $K >0$
we have $r_X(\mbt) %\rho(\mbt)^{2-Q}
= \sum_{|\smbs| < K} a(\mbs) a(\mbt + \mbs) +  \sum_{|\mbt + \smbs| < K} a(\mbs) a(\mbt + \mbs)
+ \sum_{|\smbs| \ge K, |\mbt + \smbs | \ge  K  } a(\mbs) a(\mbt + \mbs)  =: \sum_{i=1}^3  T_i(\mbt)$.
Using \eqref{abounds} and $Q > 1 $ we obtain that for any $K>0$ fixed
$|T_i(\mbt)| \le C K^3 \rho(\mbt)^{-1} = o(\rho(\mbt)^{-(2-Q)}), i=1,2 $ as $|\mbt| \to \infty $.
On the other hand, since $\liminf_{|\mbt | \to \infty}  \rho(\mbt) a(\mbt) \ge C \liminf_{|\mbt | \to \infty} g(\mbt) \ge C >0$,
see \eqref{Cpm},  \eqref{aform} so for $K >0$ large enough using \eqref{rhobdd} we infer that
\begin{eqnarray} \label{rTbounds}
T_i(\mbt)
&\ge&C\sum_{|\smbs | \ge K, |\smbt + \smbs | \ge  K  } \frac{1}{\rho (\mbs) \rho (\mbt + \mbs) } \
\ge\ C\int_{|\smbs | \ge K-3, |\smbt + \smbs | \ge  K-3  } \frac{\d \mbs }{\rho (\mbs) \rho (\mbt + \mbs) } \\
&\le&C \int_{\R^3 } \frac{\d \mbs }{\rho (\mbs) \rho (\mbt + \mbs) } + O(\rho(\mbt)^{-1}) \
= \ C(\rho^{-1} \star \rho^{-1})(\mbt) + O(\rho(\mbt)^{-1}), \nn
\end{eqnarray}
proving the lower bound in \eqref{rbounds} by \eqref{rhoQ} and \eqref{Lbounded}. The proof of the upper bound
in \eqref{rbounds} follows similarly. This proves \eqref{rbounds}.

\smallskip

Let us prove \eqref{R1cov}. Using \eqref{rbounds}, \eqref{rhobdd} and \eqref{Qh2}
we have $\sum_{\smbt \in \Z^3} |r_X(\mbt)| \ge C_2 \int_{|\smbt | > K} \rho(\mbt)^{-(2-Q)} \d \mbt
$ $=\infty $ since $2-Q < 1$.   
Next,  using \eqref{rbounds}
\begin{eqnarray*}
\sum_{(t_2,t_3)\in \Z^2} |r_X(0,t_2,t_3)| &\le&C + C\int_{[1,\infty)^2} (t_2^{q_2} + t_3^{q_3})^{-(2-Q)} \d t_2 \d t_3 \\
&\le&C +  C \int_1^\infty t_2^{-q_2(2-Q -1/q_3)} \d t_2 \int_0^\infty (1+ t_3^{q_3})^{-(2-Q)} \d t_3  < \infty
\end{eqnarray*}
follows since  $q_3(2-Q) >1 $ and
$q_2(2-Q-1/q_3) > 1 $ is equivalent to $\frac{1}{2q_1} + \sum_{j=2}^3 \frac{1}{q_j}  < 1$.
The proof of \eqref{R2cov} and \eqref{R3cov} follows similarly.
Proposition \ref{propcovZ3} is proved.
 \hfill $\Box$

\section{Limiting Gaussian random fields}

In this subsec. we define scaling limits $V^X_{\small \mbgamma}$ for $\mbgamma = (\gamma_1,\gamma_2,\gamma_3) $
satisfying \eqref{special}  (the shaded region in Fig. 1).  
The above mentioned limits are generally different
in Regions I - III of parameters $q_j,  j=1,2,3$ determined by inequalities  \eqref{R1}-\eqref{R3}.
 In some cases these limits are particularly simple and agree
with a Fractional Brownian Sheeet (FBS) with special values of Hurst parameters.

Recall that a FBS $B_{H_1,H_2,H_3}$ with parameters  $0 < H_i \le 1, i=1,2,3$ is a
Gaussian process on $\bar \R^3_+$ with zero mean and covariance
\begin{eqnarray}\label{covB}
\E \big[B_{H_1,H_2,H_3}(\mbx) B_{H_1,H_2,H_3}(\mby)\big]
&=&(1/8) \prod_{j=1}^3 (x_j^{2H_j} + y_j^{2H_j} - |x_j-y_j|^{2H_j}),
\end{eqnarray}
$\mbx = (x_1,x_2,x_3), \ \mby = (y_1,y_2,y_3) \in \bar \R^3_+ $ which is a product of the
covariances of a standard FBM with one-dimensional time parameter. Properties of FBS are discussed in \cite{aya2002}.

Let us introduce some terminology extending the terminology in \cite{ps2015}, \cite{ps2016}. If
$\ell \subset \R^3 $ is a line and  $ \mathfrak{p} \subset \R^3 $ is a plane which are orthogonal to each other,
we write $\ell \bot \mathfrak{p}$.
Also write $\mby =  (y_1,y_2,y_3) \prec \mbx =  (x_1,x_2,x_3)$ if  $y_i < x_i, i=1,2,3 $.
%from dimension $d=2$ to dimension $d = 3$.
A {\it rectangle} is a set $K(\mby, \mbx)  = \prod_{i=1}^3 (y_i, x_i] \subset \bar \R^3,
\mby \prec  \mbx,  \mbx, \mby  \in \R^3$. We say that two rectangles $ K= K(\mby, \mbx),  K' = K(\mby', \mbx'),
\mby  \prec \mbx, \mby'  \prec \mbx' $ are {\it separated by  plane}  $ \mathfrak{p} \subset \R^3 $ if
$K $ and $K'$ lie on different sides of $ \mathfrak{p}$.
By  {\it rectangular increment} of RF $V = \{ V(\mbx), \mbx \in \bar \R^3_+ \} $ on rectangle
$K(\mby, \mbx) $
we mean the (triple) difference
\begin{eqnarray*}
V(K(\mby, \mbx))&:=&V (x_1,x_2,x_3) - V (x_1,x_2,y_3) - V (x_1,y_2,x_3) - V (y_1,x_2,x_3) \\
&+&V (x_1,y_2,y_3) + V (y_1,x_2,y_3) + V (y_1,y_2,x_3) - V (y_1,y_2,y_3).
\end{eqnarray*}
We say that RF $V = \{ V(\mbx), \mbx \in \bar \R^3_+ \} $ has {\it stationary rectangular increments}  if for any
$y \in \bar \R^3_+ $,
$\{ V(K(\mby, \mbx)), \mby  \prec \mbx \} \eqfdd \{ V(K({\bf \small 0}, \mbx - \mby)), \mby  \prec \mbx \} $.

\begin{@definition}\label{def1}
Let  $V = \{ V(\mbx), \mbx \in \bar \R^3_+ \} $ be a RF with stationary rectangular increments and
$\ell \subset \R^3 $ be a line intersecting $\bar \R^3_+ $. We say that  RF $V$ has:

\smallskip

\noi (i) {\it independent rectangular increments
in direction $\ell $ } if for any orthogonal plane $\mathfrak{p} \bot \ell$
and any two rectangles $K, K' \subset \R^2_+$ separated by $\mathfrak{p}$,
the increments
$V(K)$ and $ V(K')$
are independent;

\smallskip

\noi (ii) {\it invariant rectangular increments
in direction $\ell $ } if  $V(K) = V(K')$
for any two rectangles  $K, K' \subset \R^3_+$ such that
$K' = \mbx + K$ for some $ \mbx \in \ell $.

\end{@definition}

It follows from Gaussianity and the covariance of FBS that  for $H_1 = 1/2$,
$B_{1/2,H_2,H_3}(\mbx)$ has independent rectangular increments in the direction of the coordinate axis $x_1 $ and, for $H_1 = 1$, \,
$B_{1,H_2,H_3}(x_1,x_2,x_3) = x_1 B_{H_2,H_3}(x_2,x_3) $
is a random line in $x_1$ having invariant rectangular increments in the same direction. The
case when $H_i, i=2,3 $ equal $1/2$ or 1 is analogous. Particularly,
$B_{1/2,1,H_3}$ has independent increments in direction $x_1$ and invariant increments in direction $x_2$. Except
for FBS, there are other Gaussian RFs which also enjoy the properties of increments in Definition \ref{def1}. These
RFs appear in the scaling limits of the linear RF $X$ in \eqref{Xlin} and are defined below.

Let $W(\d \mbu)$ be real-valued Gaussian white noise on $\R^3$, that it, a random process defined on Borel sets
$A \subset \R^3$ of finite Lebesgue measure  $leb(A) = \int_A \d \mbu <  \infty$ such that $W(A)$ has a Gaussian distribution with mean zero and variance
$leb(A)$ and $\E [W(A) W(B)] = leb(A\cap  B)$ for  any Borel sets $A, B \subset \R^3$ of finite Lebesgue measure.
The stochastic integral $I(g) = \int_{\R^3} g(\mbu) W(\d \mbu)$ is well-defined for any $g \in L^2(\R^3) $ and
has a Gaussian distribution $I(g) \sim N(0, \|g\|^2)$, where $\|g\|^2 = \int_{\R^3} g(\mbu)^2 \d \mbu $.
Consider the following RFs defined as stochastic integrals w.r.t. $W$:

\medskip

\noi
\begin{eqnarray}
{\cal Y}_1(\mbx)&:=&\int_{\R^3}   W(\d \mbu)
\int_{\R^3}
\frac{\1(0< u_j < x_j, j=2,3,
0< t_1 < x_1) \d \mbt}{ (c_1|t_1-u_1|^{\frac{q_1}{\nu}} + \sum_{j=2}^3
c_j|t_j|^{\frac{q_j}{\nu}})^{\nu}}, \hskip1cm \text{(Region I)}\label{R1def} \\
%B_{1, {\cal H}_2, 1/2 }(x_1,x_2,x_3)
{\cal Y}_2(\mbx)
&:=&x_1 \int_{\R^3} W(\d \mbu)
\int_{\R^2}
\frac{ \1(0< u_3 < x_3, 0< t_2 < x_2) \d t_2 \d t_3}{ (c_1|u_1|^{\frac{q_1}{\nu}} + c_2|t_2-u_2|^{\frac{q_2}{\nu}} +
c_3|t_3|^{\frac{q_3}{\nu}})^{\nu}}, \hskip.7cm \text{(Region II)}
\label{R2def} \\
%B_{1,1,{\cal H}_3}&=&x_1 x_2  \int_{\R^3} Z(\d \mbu)
{\cal Y}_3(\mbx)
&:=&x_1x_2\int_{\R^3} W(\d \mbu)
\int_\R
\frac{\1(0<t < x_3)\d t}{ (\sum_{j=1}^2 c_j|u_j|^{\frac{q_j}{\nu}} +
c_3|t - u_3|^{\frac{q_3}{\nu}})^{\nu}}, \hskip1.5cm  \text{(Region III)}\label{R3def} \\
{\cal Y}_{12}(\mbx)
&:=&\int_{\R^3} W(\d \mbu)
\int_{\R^3}
\frac{\1(0< t_j < x_j, j=1,2, 0< u_3 < x_3) \d \mbt}{
(\sum_{i=1}^2 c_j|t_j-u_j|^{\frac{q_j}{\nu}} +
c_3|t_3|^{\frac{q_3}{\nu}})^\nu}, \hskip.8cm  \text{(Regions I\&II)} \label{R4def}  \\
{\cal Y}_{23}(\mbx)
&:=&x_1\int_{\R^3} W(\d \mbu)
\int_{\R^2}
\frac{\1(0< t_j < x_j, j=2,3) \d t_2 \d t_3}{
(c_1|u_1|^{\frac{q_1}{\nu}} + \sum_{j=2}^3 c_j|t_j-u_j|^{\frac{q_j}{\nu}})^\nu}, \hskip1.2cm  \text{(Regions II\&III)}
\label{R5def}\\
{\cal Y}_0(\mbx)
&:=&\int_{\R^3} W(\d \mbu)
\int_{\R^3}
\frac{\1(0< t_j < x_j, j=1,2,3) \d \mbt}{
(\sum_{j=1}^3 c_j|t_j-u_j|^{\frac{q_j}{\nu}})^\nu}. \hskip2.2cm  \text{(Regions I\&II\&III)} \label{R6def}
\end{eqnarray}
We also write ${\cal Y}_1(\mbx) \equiv {\cal Y}_1(\mbx;  \mbq, \mbc), \cdots,
{\cal Y}_0(\mbx) \equiv {\cal Y}_0(\mbx;  \mbq, \mbc)$ to emphasize the dependence of these RFs on vector parameters
$\mbq = (q_1,q_2,q_3) \in \R^3_+ $ and
$\mbc = (c_1,c_2,c_3) \in \R^3_+$.

\begin{thm}\label{exist} (i) The Gaussian RFs in \eqref{R1def}-\eqref{R6def} depending on vector parameters $\mbq = (q_1,q_2,q_3) \in \R^3_+ $ and
$\mbc = (c_1,c_2,c_3) \in \R^3_+$
are well-defined in the indicated parameter regions given in \eqref{R1}-\eqref{R3}. More precisely,

\bigskip

\noi (i1) \ ${\cal Y}_1$ in \eqref{R1def} is well-defined for \ \
$ \frac{1}{2q_1} + \sum_{i=2}^3 \frac{1}{q_i}  < 1 <  \sum_{i=1}^3 \frac{1}{q_i}.$

\bigskip

\noi (i2) \ ${\cal Y}_2$ in \eqref{R2def} is well-defined for \ \
$\sum_{i=1}^2 \frac{1}{2q_i} + \frac{1}{q_3} < 1 < \frac{1}{2q_1} + \sum_{i=2}^3 \frac{1}{q_i}.$

\bigskip

\noi (i3) \ ${\cal Y}_3$ in \eqref{R3def} is well-defined for \ \
$\sum_{i=1}^3 \frac{1}{2q_i} < 1 < \sum_{i=1}^2 \frac{1}{2q_i} + \frac{1}{q_3}.$

\bigskip

\noi (i4) \ ${\cal Y}_{12}$  in \eqref{R4def}  is well-defined for \ \
$\sum_{i=1}^2 \frac{1}{2q_i} + \frac{1}{q_3} < 1 <  \sum_{i=1}^3 \frac{1}{q_i}.$

\bigskip

\noi (i5) \ ${\cal Y}_{23}$  in \eqref{R5def}  is well-defined for \ \
$\sum_{i=1}^3 \frac{1}{2q_i} < 1 <  \frac{1}{2q_1} + \sum_{i=2}^3 \frac{1}{q_i}.$

\bigskip

\noi (i6) \ ${\cal Y}_0$  in \eqref{R6def}  is well-defined for \ \
$\sum_{i=1}^3 \frac{1}{2q_i} < 1 <  \sum_{i=1}^3 \frac{1}{q_i}.$

\bigskip

\noi (ii) RFs in \eqref{R1def}-\eqref{R6def} have stationary rectangular increments and satisfy the self-similarity properties:
\begin{eqnarray*}
{\cal Y}_1(\lambda_1 x_1, \lambda_2 x_2, \lambda_3 x_3)&\eqfdd&\lambda_1^{{\cal H}_1} \lambda_2^{1/2} \lambda_3^{1/2}
{\cal Y}_1(x_1,x_2, x_3),
\hskip1cm \forall \lambda_i >0, \ i=1,2,3, \\
{\cal Y}_2(\lambda_1 x_1, \lambda_2 x_2, \lambda_3 x_3)&\eqfdd&\lambda_1 \lambda_2^{{\cal H}_2}\lambda_3^{1/2}
{\cal Y}_2(x_1,x_2, x_3), \hskip1.4cm \forall \lambda_i >0, \ i=1,2,3, \\
{\cal Y}_3(\lambda_1 x_1, \lambda_2 x_2, \lambda_3 x_3)&\eqfdd&\lambda_1 \lambda_2 \lambda_3^{{\cal H}_3}
{\cal Y}_3(x_1,x_2, x_3),
\hskip1.7cm \forall \lambda_i >0, \ i=1,2,3, \\
{\cal Y}_{12}(\lambda^{1/q_1} x_1, \lambda^{1/q_2} x_2, \mu x_3)&\eqfdd&
\lambda^{{\cal H}_{12}}
\mu^{1/2}
{\cal Y}_{12}(x_1,x_2, x_3),
\hskip1.5cm \forall \lambda >0, \ \mu >0, \\
{\cal Y}_{23}(\lambda x_1, \mu^{1/q_2} x_2, \mu^{1/q_3} x_3)&\eqfdd&
\lambda \mu^{{\cal H}_{23}}
{\cal Y}_{23}(x_1,x_2, x_3),
\hskip2cm \forall \lambda >0, \ \mu >0, \\
{\cal Y}_0(\lambda^{1/q_i} x_i, i=1,2,3)&\eqfdd&
\lambda^{{\cal H}_0}{\cal Y}_0(x_1,x_2, x_3),
\hskip2.5cm \forall \lambda >0,
\end{eqnarray*}
where
\begin{eqnarray}\label{calH}
&{\cal H}_1\ :=\ \frac{3}{2} - q_1 (1- \frac{1}{q_2} - \frac{1}{q_3}), \qquad
{\cal H}_2\ :=\ \frac{3}{2} - q_2(1- \frac{1}{2q_1} - \frac{1}{q_3}), \qquad
{\cal H}_3\ :=\  \frac{3}{2} - q_3(1- \frac{1}{2q_1} - \frac{1}{2q_2}), \\
&{\cal H}_{12}\ :=\ \frac{3}{2q_1} + \frac{3}{2q_2} + \frac{1}{q_3} -1, \qquad
{\cal H}_{23}\ :=\  \frac{1}{2q_1} + \frac{3}{2q_2} + \frac{3}{2q_3} -1, \qquad
{\cal H}_0\ :=\  \sum_{i=1}^3 \frac{3}{2q_i} - 1. \nn
\end{eqnarray}

\medskip

\noi (iii)\, RFs ${\cal Y}_i, i=1,2,3$ agree, up to multiplicative constants $\sigma_i := \E^{1/2} [{\cal Y}^2_i (1,1,1; \mbq, \mbc)]$,
with  FBS %$B_{{\cal H}_1,{\cal H}_2, {\cal H}_3}$
having  two its parameters  %$0 < {\cal H}_i \le 1, i=1,2,3$
equal to either $1/2$ or 1. Namely,
\begin{eqnarray}\label{YFBS}
{\cal  Y}_1&\eqfdd&\sigma_1 B_{{\cal H}_1,1/2,1/2}, \qquad %\mbox{${\cal H}_1  =\frac{3}{2} - q_1 (1- \frac{1}{q_2} - \frac{1}{q_3}) = \frac{1+ H_1}{2}$}, \\
{\cal  Y}_2\ \eqfdd \ \sigma_2 B_{1, {\cal H}_2, 1/2 }, \qquad %\mbox{${\cal H}_2 =   \frac{3}{2} - q_2(1- \frac{1}{2q_1} - \frac{1}{q_3}) =
%\frac{1}{2} +  \frac{H_2(H_1-1)}{2H_1}$}, \\
{\cal  Y}_3\ \eqfdd\ \sigma_3 B_{1,1, {\cal  H}_3}, %\hskip1.5cm   \mbox{${\cal H}_3 := \frac{1}{2} +  \frac{H_3}{2H_1H_2}(H_1H_2 - H_1 - H_2) $}
\end{eqnarray}
where ${\cal H}_i, i=1,2,3 $ are defined in \eqref{calH}.

\end{thm}

\medskip

\noi {\it Proof.}  (i) In view of inequalities \eqref{Cpm} and the form of the integrands,
it suffices prove the existence of the stochastic integrals
for $c_1=c_2 = c_3 = \nu =  1 $ and $x_1 = x_2 = x_3 = 1$.

\smallskip

\noi (i1) It suffices to prove
\begin{eqnarray*}
I&:=&\int_{\R} \d u
\Big( \int_0^{1} \int_{\R} \int_{\R}
\frac{\d t_1 \d t_2 \d t_3}{|t_1-u|^{q_1} + |t_2|^{q_2} +
|t_3|^{q_3}} \Big)^2 \ < \ \infty.
\end{eqnarray*}
Split $I = I_1 + I_2$, where $I_1 := \int_{|u| < 2} \cdots, \ I_2 := \int_{|u| > 2} \cdots$. Then
%\begin{eqnarray*}
$ I_1\le C
\big(\int_0^{1} \int_0^\infty \int_0^\infty (t_1^{q_1} + t_2^{q_2} +
t_3^{q_3})^{-1} \d t_1 \d t_2 \d t_3 \big)^2.
$ %\end{eqnarray*}
Using twice the second inequality in \eqref{elem}, we obtain
%We have by change of variables $t_2 \to t^{q_1/q_2} t_2, t_3 \to t^{q_1/q_3} t_3$ that
\begin{eqnarray*}
&J(t) := \int_0^\infty \int_0^\infty (t_1^{q_1} + t_2^{q_2} +
t_3^{q_3})^{-1} \d t_2 \d t_3  \ \le  C \int_0^\infty (t_1^{q_1} + t_2^{q_2})^{(1/q_3)-1} \d t_2 \le C
 \  t_1^{(q_1/q_2) + (q_1/q_3) - q_1},
\end{eqnarray*}
and hence $I_1 %\le C \big(\int_0^1 t^{-q_1(1 - \frac{1}{q_2} - \frac{1}{q_3})} \d t \big)^2
<  \infty $
since $q_1(1 - \frac{1}{q_2} - \frac{1}{q_3}) < 1 $. % is equivalent to $q_1q_2 q_3 < q_1 q_2 + q_1 q_3 + q_2 q_3 $.
Similarly,
%\begin{eqnarray*}
$I_2 \le C\int_1^\infty \d u
\big(\int_0^\infty \int_0^\infty (u^{q_1} + t_2^{q_2} +
t_3^{q_3})^{-1} \d t_2 \d t_3 \big)^2 \ = \ C \int_1^\infty
J^2(u) \d u  < \infty $
%\le\  C\int_1^\infty u^{-2q_1(1 - \frac{1}{q_2} - \frac{1}{q_3})}\d u  \ < \ \infty
%\end{eqnarray*}
since $2q_1(1 - \frac{1}{q_2} - \frac{1}{q_3}) > 1 $. This proves that ${\cal Y}_1$ in \eqref{R1def} is well-defined.

\smallskip

\noi (i2) It suffices to prove $I:= \int_{\R^2} \d u_1 \d u_2 \big(  \int_0^{1} I(u_1, t_2-u_2) \d t_2\big)^2 < \infty $, where 
\begin{eqnarray}\label{Iuv}
&I(u,v) := \int_{\R}
(|u|^{q_1} + |v|^{q_2} + |t|^{q_3})^{-1} \d t   \ \le \  C(|u|^{q_1} + |v|^{q_2})^{\frac{1}{q_3}-1},
\end{eqnarray}
see  \eqref{elem}.
%where $K := \int_{\R} (1 + |t|^{q_3})^{-1} \d t < \infty $ since $q_3 > 1$.
Therefore, $I \le  C \int_{\R^2} \d u_1 \d u_2 \big(\int_0^1 (|u_1|^{q_1} + |t_2 - u_2|^{q_2})^{\frac{1}{q_3}-1} \d t_2 \big)^2
=  C \big( \int_{\R} \d u_1 \int_{|u_2|< 2} \d u_2 (\cdots)^2 + $   $ \int_{\R} \d u_1 \int_{|u_2|> 2} \d u_2 (\cdots)^2 \big)
=: C(I_1 + I_2)$, where $ I_1 \le  C \int_{0}^\infty F(u_1)^2 \d u_1 $ and
\begin{eqnarray*}
&F(u_1) := \int_0^1 (u_1^{q_1} + t_2^{q_2})^{\frac{1}{q_3}-1} \d t_2
\le C \begin{cases} u_1^{-q_1(1- \frac{1}{q_2} - \frac{1}{q_3})}, &1> \frac{1}{q_2} - \frac{1}{q_3}, \\
1, &1< \frac{1}{q_2} - \frac{1}{q_3}, \\
|\log u_1|,  &1= \frac{1}{q_2} - \frac{1}{q_3}
\end{cases}
\end{eqnarray*}
according to \eqref{elem}, implying $I_1 \le  C \int_{0}^\infty F(u_1)^2 \d u_1 < \infty $.  
We also have
$I_2 =\int_1^\infty  \d u_2  \int_0^\infty (u_1^{q_1} + u_2^{q_2})^{2(\frac{1}{q_3}-1)}  \d u_1
\le  \int_1^\infty  u_2^{- 2q_2(1- \frac{1}{q_3}) - \frac{q_2}{q_1}  } \d u_2  <  \infty $
%\end{eqnarray*}
since $2q_2(1- \frac{1}{q_3}) - \frac{q_2}{q_1}> 1$. This proves that ${\cal Y}_2$ in \eqref{R2def} is well-defined.

\smallskip

\noi (i3) It suffices to prove
%\begin{eqnarray*}
$I := \int_{\R^3} \d u_1 \d u_2 \d u_3
\big( \int_0^{1}
(u_1^{q_1} + |u_2|^{q_2} +
|t-u_3|^{q_3})^{-1} \d t \big)^2
%\end{eqnarray*}
 = \int_{\R^2} \d u_1 \d u_2 \int_{|u_3|<2} \d u_3 (\dots)^2 +  \int_{\R^2} \d u_1 \d u_2  \int_{|u_3|>2} \d u_3 (\dots)^2
=: I_1 + I_2 < \infty.$   
We have $I_1 \le C\int_{\R^2_+} G(u_1,u_2)^2 \d u_1 \d u_2,$ where $
G(u_1,u_2) :=  \int_0^1 (u_1^{q_1} + u_2^{q_2} +t^{q_3})^{-1} \d t$ can be estimated by \eqref{elem}.
Using \eqref{elem}
we obtain $\int_{ u_1^{q_1} + u_2^{q_2} > 1}  G(u_1,u_2)^2 \d u_1 \d u_2
\le C \int_{\R^2_+} (u_1^{q_1} + u_2^{q_2} + 1)^{-2} \d u_1 \d u_2
\le C \int_0^\infty (u_1^{q_1} + 1)^{\frac{1}{q_2}-2} \d u_1 < \infty $
since $(2-\frac{1}{q_2})q_1 > 1 $. Using the same inequality, in the case  $2q_2(1-\frac{1}{q_2}) > 1 $
we conclude that
$\int_{ u_1^{q_1} + u_2^{q_2} \le 1} G(u_1,u_2)^2 \d u_1 \d u_2 \le
\int_{[0,1]^2} G(u_1,u_2)^2 \d u_1 \d u_2 \le $  $ C \int_0^1 u_1^{-2q_1(1- \frac{1}{2q_2} - \frac{1}{q_3})} \d u_1 $  $ < \infty $
since  $2q_1(1- \frac{1}{2q_2} - \frac{1}{q_3}) < 1 $; in the case  $2q_2(1-\frac{1}{q_2}) > 1 $ the convergence
of $\int_{[0,1]^2} G(u_1,u_2)^2 \d u_1 \d u_2 $ follows trivially from \eqref{elem}. This proves
$I_1 < \infty $.   
Finally, $I_2 \le C\int_1^\infty  \d u_3 \int_{\R^2_+} (u_1^{q_1} + u_2^{q_2} +
u_3^{q_3})^{-2} \d u_1 \d u_2
\le C\int_1^\infty  \d u_3  \int_0^\infty
(u_1^{q_1} + u_3^{q_3})^{-(2- \frac{1}{q_2})} \d u
\le C\int_1^\infty
u_3^{-q_3(2- \frac{1}{q_1} -\frac{1}{q_2})} \d z   <  \infty $ 
since $q_3(2- \frac{1}{q_1} -\frac{1}{q_2})>1$.  
This proves that ${\cal Y}_3 $  in \eqref{R3def} is well-defined.

\medskip

\noi (i4) It suffices to show
%\begin{eqnarray*}
$I :=\int_{\R^2} \d u_1 \d u_2  \big( \int_{[0,1]^2 \times \R}
(|t_1-u_1|^{q_1} + |t_2-u_2|^{q_2} +
|u_3|^{q_3})^{-1} \d t_1 \d t_2 \d u_3 \big)^2  <  \infty. $
%\end{eqnarray*}
Using \eqref{elem},
%\begin{eqnarray*}
$I\le  C\int_{\R^2} \d u_1 \d u_2 \big( \int_{[0,1]^2}
(|t_1-u_1|^{q_1} + |t_2-u_2|^{q_2})^{\frac{1}{q_3} -1}
\d t_1 \d t_2  \big)^2 =  C\sum_{k=1}^4 I_k,$,
%\Big( \int_{[0,1]^2}
%\frac{\d t \d s }{(|t-u|^{q_1} + |s-v|^{q_2})^{1- \frac{1}{q_3}}} \Big)^2  =  C\sum_{k=1}^4 I_k,
%\end{eqnarray*}
where $I_1 := \int_{|u_1|\le 2, |u_2|\le 2} \d u_1 \d u_2 \, (\cdots )^2,$ \
$I_2 := \int_{|u_1|\le 2, |u_2| > 2} \d u_1 \d u_2 \, (\cdots)^2,$ \
$I_3 := \int_{|u_1| > 2, |u_2|\le 2} \d u_1 \d u_2 \, (\cdots)^2, $ \  and
$I_4 :=  $   $\int_{|u_1| > 2, |u_2| > 2} \d u_1 \d u_2 \, (\cdots)^2$. Here similarly to
the proof of $\int_{[0,1]^2} G(u_1,u_2)^2 \d u_1 \d u_2 < \infty $ in (i3), in the case $q_2(1-\frac{1}{q_3})> 1 $ we obtain that
%\begin{eqnarray*}
$I_1\le C \big(\int_{[0,1]^2}  (t_1^{q_1} + t_2^{q_2})^{\frac{1}{q_3}-1} \d t_1 \d t_2 \big)^2 \le
%\frac{ t \d s }|t-u|^{q_1} + |s-v|^{q_2})^{1- \frac{1}{q_3}}} \le C
%\begin{cases}1,  &q_2(1-\frac{1}{q_3})< 1, \\
\big(\int_0^1 t_1^{-q_1(1- \frac{1}{q_2} - \frac{1}{q_3})} \d t_1 \big)^2 < \infty $
%\end{cases}
%\end{eqnarray*}
since $1 < \sum_{i=1}^3 \frac{1}{q_i}$, and in the case $q_2(1-\frac{1}{q_3})\le 1 $ the same result
$I_1 < \infty $ follows even easier. Next, by \eqref{elem} 
$I_2\le C\int_1^\infty \d u_2 \big( \int_0^{1} (t_1^{q_1} + u_2^{q_2})^{\frac{1}{q_3}-1} \d t_1 \big)^2  \le  C \int_1^\infty
u_2^{-2q_2(1- \frac{1}{q_3})} \d u_2 < \infty $ 
since $2q_2(1- \frac{1}{q_3})> 1$. Similarly,
$I_3\le C \int_1^\infty u_1^{-2q_1(1- \frac{1}{q_3})} \d u_1 < \infty $. Finally,  
$I_4 \le C\int_{[1,\infty)^2} (u_1^{q_1} + u_2^{q_2})^{-2(1- \frac{1}{q_3})} \d u_1 \d u_2  \le  C \int_1^\infty
u_1^{-2q_1(1- \frac{1}{2q_2}  - \frac{1}{q_3})} \d u_1 < \infty $ 
as $1>  \frac{1}{2q_1} + \frac{1}{2q_2}  + \frac{1}{q_3}$. This proves
that ${\cal Y}_{12}$ in \eqref{R4def} is well-defined.

\medskip

\noi (i5) It suffices to show
$I := \int_{\R^3} \d u_1 \d u_2 \d u_3
\big( \int_{[0,1]^2}
(|u_1|^{q_1} + |t_2- u_2|^{q_2} +
|t_3- u_3|^{q_3})^{-1} \d t_2 \d t_3 \big)^2 \ < \ \infty. $   
Split $I= \sum_{j=1}^4 I_j$, where $I_1 := \int_\R \d u_1 \int_{|u_2|\le 2, |u_3| \le 2} \d u_2 \d u_3 (\cdots)^2, \
I_2 := \int_\R \d u_1 \int_{|u_2| \le 2, |u_3| > 2} \d u_2 \d u_3  (\cdots)^2, \
I_3 := $   $ \int_\R  \d u_1 $    $\int_{|u_2| > 2, |u_3| \le 2} \d u_2 \d u_3 (\cdots)^2, \
I_4 := \int_\R \d u_1 \int_{|u_2|>2, |u_3|>2} \d u_2 \d u_3  (\cdots)^2 $. Then using \eqref{elem} in the case
$q_3 > 1$ we obtain
%\begin{eqnarray*}
$I_1 \le C\int_0^\infty \d u_1
\big( \int_{[0,1]^2}
(u_1^{q_1} + \sum_{i=2}^3 t_i^{q_i})^{-1} \d t_2 \d t_3 \big)^2
\le C \int_1^\infty u_1^{-2q_1} \d u_1 + C \int_0^1 \d u_1
\big(\int_0^{1} (u^{q_1} +  t_2^{q_2})^{-(1- \frac{1}{q_3})} \d t_2 \big)^2
$   $\le C  + C \int_0^1 u_1^{-2q_1 (1  - \frac{1}{q_2} - \frac{1}{q_3})} \d u_1  < \infty $  
since $2q_1 (1  -  \frac{1}{q_3}) + 2\frac{ q_1}{q_2} < 1 $; when $q_3 \le 1 $  the convergence $I_{1} < \infty $
follows easily.
%is equivalent to $q_1q_2 q_3 < q_1 q_2 + q_1 q_3 +  \frac{q_2q_3}{2}$.
%(In the case when  $q_3 > 1> \frac{1}{q_2} + \frac{1}{q_3} $ is violated, $I_{11} < \infty $ follows from  \eqref{kappabdd}, too.)
Hence, $I_1 < \infty$.  Next, 
$I_3  \le C  \int_0^\infty \d u \int_1^\infty (u^{q_1} + v^{q_2})^{2} \d v \le C \int_0^\infty v^{-2q_2(1- \frac{1}{2q_1})} \d v < \infty $ since
$1> \sum_{i=1}^2 \frac{1}{2q_i}$ and 
$I_2 < \infty $ follows analogously. Finally, due to $1> \sum_{i=1}^3 \frac{1}{2q_i}$, relations  
$I_4 \le C \int_{\R_+ \times [1, \infty)^2}  (u_1^{q_1} + u_2^{q_2} + u_3^{q_3})^{-2} \d u_1 \d u_2 \d u_3
\le C \int_{[1,\infty)^2} (u_1^{q_1} + u_2^{q_2})^{-(2- \frac{1}{q_3})} \d u_1 \d u_2
%\big(\int_0^{1} (u^{q_1} +  t^{q_2})^{-(1- \frac{1}{q_3})} \d t \big)^2   \
%\le\  C\int_0^1 u^{-2q_1 (1  - \frac{1}{q_2} - \frac{1}{q_3})} \d u \big(\int_0^{1/u^{q_1/q_2}} (1+ t^{q_2})^{-(1- \frac{1}{q_3})} \d t \big)^2
 <  \infty $ follow
similarly as in the proof of $I_4 < \infty $ in (i4). This proves
that ${\cal Y}_{23}$ in \eqref{R5def} is well-defined.

\medskip

\noi (i6) It suffices to show
\begin{eqnarray*}
&I \ := \ \int_{\R^3} \d u_1 \d u_2 \d u_3
\Big( \int_{[0,1]^3}
\frac{\d t_1 \d t_2 \d t_3}{\sum_{i=1}^3 |t_i-u_i|^{q_i}} \Big)^2 \ < \ \infty.
\end{eqnarray*}
Split $I= \sum_{k=1}^8 I_k$ into the sum of 8 integrals according to whether
$|u_i| \le 2 $ or $|u_i| > 2, i=1,2,3 $.  In the case $1> \frac{1}{q_2} + \frac{1}{q_3} $ using
\eqref{elem}  we obtain
\begin{eqnarray*}
I_1\le&C\int_{[0,1]^3}
(\sum_{i=1}^3 t_i^{q_i})^{-1} \d t_1 \d t_2 \d t_3
 \ \le \ C \int_0^1 \int_0^1  \big(\sum_{i=1}^2 t_i^{q_i} \big)^{-(1- \frac{1}{q_3})}
\d t_1 \d t_2 \ \le \ C \int_0^1 t^{-q_1(1- \frac{1}{q_2} - \frac{1}{q_3})} \d t \ < \ \infty;
\end{eqnarray*}
for $1\ge \frac{1}{q_2} + \frac{1}{q_3} $ relation $I_1 < \infty $ follows easily. The remaining
integrals can be easily evaluated, e.g,
\begin{eqnarray*}
&&\mbox{$I_2\le C\int_1^\infty \d u_1 \big(\int_0^1 \int_0^1 \frac{\d t_2 \d t_3}
{u_1^{q_1} + \sum_{i=2}^3 t_i^{q_i}} \big)^2 \le C \int_1^\infty u_1^{-2q_1} \d u_1 < \infty,$} \\
&&\mbox{$I_3\le C\int_1^\infty \int_1^\infty \d u_1 \d u_2 \big(\int_0^1 \frac{\d t_3}
{\sum_{i=1}^2 u_i^{q_i} + t_3^{q_3}} \big)^2 \le C \int_1^\infty \int_1^\infty \frac{\d u_1 \d u_2}
{\big(\sum_{i=1}^2 u_i^{q_i}\big)^2} \ \le \ \int_1^\infty  u_1^{-q_1(2-\frac{1}{q_2})} \d u_1  < \infty,$}\\
&&\mbox{$I_4\le C\int_1^\infty \int_1^\infty \int_1^\infty \frac{\d u_1 \d u_2 \d u_3}
{(\sum_{i=1}^3 u_i^{q_i})^2} \le C\int_1^\infty \int_1^\infty  \frac{\d u_1 \d u_2 }
{(\sum_{i=1}^2 u_i^{q_i})^{2 - \frac{1}{q_3} } } \le C\int_1^\infty
u_1^{-q_1(2 - \frac{1}{q_2} - \frac{1}{q_3}) } \d u_1 < \infty.$}
\end{eqnarray*}
This proves
that ${\cal Y}_{0}$ in \eqref{R6def} is well-defined, thereby completing the proof of part (i).

\medskip

\noi (ii)  The self-similarity properties follow from scaling properties
$\{ W(\d \lambda_1 u_1, \d \lambda_2 u_2,\d \lambda_3 u_3) \} \eqfdd \{
(\lambda_1 \lambda_2 \lambda_3)^{1/2} $   $ W(\d u_1, \d u_2,\d u_3) \}\, (\forall \lambda_i >0, i=1,2,3) $ of the white noise
and the integrands in  \eqref{R1def}-\eqref{R6def}. For example,
\begin{eqnarray*}
&&{\cal Y}_{12}((\lambda^{1/q_1} x_1, \lambda^{1/q_2} x_2, \mu x_3;\mbq, \mbc)\\
&=&\lambda^{\sum_{i=1}^3 1/q_i}\int_{\R^3} W(\d \lambda^{1/q_1} u_1, \d \lambda^{1/q_2} u_2, \d \mu u_3)
\int_{\R^3}
\frac{\1(0< \lambda^{1/q_i} t_i < \lambda^{1/q_i}x_i, i=1,2, 0< \mu u_3 < \mu x_3)  \d \mbt}{
(\sum_{i=1}^2 c_i|\lambda^{1/q_i} t_i- \lambda^{1/q_i} u_i|^{\frac{q_i}{\rho}} +
c_3|\lambda^{1/q_3} t_3|^{\frac{q_3}{\nu}})^\nu} \\
&\eqfdd&\lambda^{{\cal H}_{12}} \mu^{1/2}{\cal Y}_{12}(x_1, x_2, x_3;\mbq, \mbc).
\end{eqnarray*}

\medskip

\noi (iii)  By Gaussianity, it suffices to show the agreement of the corresponding covariance functions.
%${\cal Y}_i(\mbx) \equiv {\cal Y}_i(\mbx;\mbq, \mba), i=1,2,3$.
Using the definition in \eqref{R1def} we have that
$\E [{\cal Y}_1(\mbx; \mbq, \mbc) {\cal Y}_1(\mby; \mbq, \mbc)]  = \int_0^{x_1} \int_0^{y_1} \theta(t-s) \d t \d s \prod_{i=2}^3 (x_i \wedge y_i), $
where
\begin{eqnarray*}
\theta(t)
&:=&\int_{\R^5}
\frac{\d u \, \d t_2 \, \d t_3 \, \d s_2 \, \d s_3}{ (c_1|u|^{\frac{q_1}{\rho}} + \sum_{i=2}^3 c_i|t_i|^{\frac{q_i}{\nu}})^\nu
(c_1|t-u|^{\frac{q_1}{\rho}} + \sum_{i=2}^3 c_i|s_i|^{\frac{q_i}{\nu}})^{\nu} } \ = \  \theta(1)|t|^{1+ 2q_1 (\frac{1}{q_2} + \frac{1}{q_3} -1)}.
\end{eqnarray*}
Hence using $3+ 2q_1 (\frac{1}{q_2} + \frac{1}{q_3} -1)  = 2{\cal H}_1$ we obtain
$$
\int_0^{x} \int_0^y \theta(t-s)  \d t \d s = (C_1/2)(x^{2{\cal H}_1} + y^{2{\cal H}_1}
- |x-y|^{2{\cal H}_1}), \quad x, y \ge 0,
$$
proving $\E [{\cal Y}_1(\mbx; \mbq, \mbc) {\cal Y}_1(\mby; \mbq, \mbc)]
= C_1 \E [B_{{\cal  H}_1, 1/2, 1/2}(\mbx)B_{{\cal  H}_1, 1/2, 1/2}(\mby)], \, \mbx, \mby \in \R^3_+,$
for some constant
$C_1 >0$. Particularly, $C_1 \E [B^2_{{\cal  H}_1, 1/2, 1/2}(1,1,1)] = C_1 = \E [{\cal Y}^2_1(1,1,1)]$, or $C_1 = \sigma_1^2$.
This proves the first relation in \eqref{YFBS} and the other two relations \eqref{YFBS} follow analogously.
Theorem \ref{exist} is proved. \hfill $\Box$

\begin{remark} The self-similarity properties in (ii) imply the following
operator scaling properties of the corresponding RFs.
%can be rewritten as follows.
For $\lambda >0, \mbgamma = (\gamma_1,\gamma_2,\gamma_3) \in \R^3_+$ denote
the diagonal $3\times 3$-matrix
$\lambda^\mbgamma = {\rm diag}( \lambda^{\gamma_i}, i=1,2,3).$ Then for any $\lambda>0$
\begin{eqnarray}
{\cal Y}_i(\lambda^{\small \mbgamma} \mbx;  \mbq, \mba)&\eqfdd&\lambda^{H_i ({\small \mbgamma, \mbq})}
{\cal Y}_i(\mbx;  \mbq, \mba),
\hskip1cm  i=1,2,3, \label{y1}\\
{\cal Y}_{ij}(\lambda^{\small \mbgamma}\mbx;  \mbq, \mba)&\eqfdd&\lambda^{H_{ij} ({\small \mbgamma, \mbq})}
{\cal Y}_{ij}(\mbx;  \mbq, \mba), \quad \gamma_i q_i = \gamma_j q_j,
\hskip.3cm  1\le i<j\le 3, j = i+1, \label{y2}\\
{\cal Y}_0(\lambda^{\small \mbgamma}\mbx;  \mbq, \mba)&\eqfdd&\lambda^{H_0 ({\small \mbgamma, \mbq})}
{\cal Y}_0(\mbx;  \mbq, \mba),\quad \gamma_1q_1 = \gamma_2 q_2 = \gamma_3 q_3, \label{y0}
\end{eqnarray}
where
\begin{eqnarray}\label{H}
&&\mbox{$H_1(\mbgamma, \mbq)\ :=\ \gamma_1 {\cal H}_1+ \frac{\gamma_2+\gamma_3}{2} \ = \ \frac{3\gamma_1+ \gamma_2 + \gamma_3}{2} +
\gamma_1 q_1 (\frac{1}{q_2} + \frac{1}{q_3} - 1)$}, \\
&&\mbox{$H_2(\mbgamma, \mbq)\ :=\
\gamma_1  + \gamma_2{\cal H}_2+ \frac{\gamma_3}{2} \ = \  \gamma_1 + \frac{3\gamma_2 + \gamma_3}{2} +
\gamma_2 q_2(\frac{1}{2q_1} + \frac{1}{q_3} - 1)$}, \nn \\
&&\mbox{$H_3(\mbgamma, \mbq)\ :=\ \gamma_1 + \gamma_2 + \gamma_3{\cal H}_3 \ = \ \gamma_1 + \gamma_2 + \frac{3\gamma_3}{2} +
\gamma_3q_3 (\frac{1}{2q_1} + \frac{1}{2q_2}- 1)$}, \nn \\
&&\mbox{$H_{12}(\mbgamma, \mbq)\ :=\ \gamma_1 q_1 {\cal H}_{12}+ \frac{\gamma_3}{2} \ = \ \frac{3(\gamma_1 + \gamma_2)+ \gamma_3}{2} +
\gamma_1 q_1 (\frac{1}{q_3}-1)$},  \nn \\
&&\mbox{$H_{23}(\mbgamma, \mbq)\ :=\ \gamma_1 + \gamma_2 q_2 {\cal H}_{23} \ = \  \gamma_1 + \frac{3\gamma_2 + 3\gamma_3}{2} +
\gamma_2 q_2(\frac{1}{2q_1} - 1)$}, \nn \\
&&\mbox{$H_0(\mbgamma, \mbq)\ :=\ \gamma_1 q_1 {\cal H}_0 \ = \ \gamma_1 q_1 (\sum_{i=1}^3 \frac{3}{2q_i} -1). $} \nn
\end{eqnarray}
See \cite{bier2007} for the definition and general properties of operator scaling RFs.
Note that  while \eqref{y1} hold for any $\mbgamma \in \R^3_+$, the self-similarity properties
in \eqref{y2}  and \eqref{y0} hold for $\mbgamma \in \R^3_+$ satisfying one and two (all)  balance conditions
in \eqref{bal}, respectively. Also note that $H_1(\mbgamma, \mbq) = H_2(\mbgamma, \mbq) = H_{12}(\mbgamma, \mbq) $ 
for $\gamma_1 q_1 = \gamma_2 q_2$,  $H_2(\mbgamma, \mbq) = H_3(\mbgamma, \mbq) = H_{23}(\mbgamma, \mbq) $
for $\gamma_2 q_2 = \gamma_3 q_3$, and that {\it all} scaling exponents in \eqref{H} coincide for 
$\gamma_1 q_1 = \gamma_2 q_2 = \gamma_3 q_3.$
\end{remark}

\begin{remark} It follows from \eqref{YFBS} that RFs ${\cal Y}_i, i=1,2,3 $ in \eqref{R1def}-\eqref{R3def} have the 
(rectangular) increment properties of Definition \ref{def1} in {\it two} directions in $\R^3$.
For instance, ${\cal Y}_1 $ has independent increments in $x_2 $ and $x_3$, while
${\cal Y}_3 $ has invariant increments in $x_1 $ and $x_2$. The RFs
${\cal Y}_{12} $ and ${\cal Y}_{23} $  have these properties in {\it one} direction,
namely, ${\cal Y}_{12} $ has independent increments in $x_3$ and
${\cal Y}_{23} $ has invariant increments in $x_3$. These facts follow from the representations
in \eqref{R4def} and \eqref{R5def} and the independent increment property of the white noise
$W(\d \mbu)$.
They are closely related to the number of  balance conditions satisfied by $\mbgamma$'s
as shown in the following sec.

\end{remark}

\section{The main result}

In this sec. we formulate our main result about partial sums limits in \eqref{Sn} of the linear RF
$X$ in \eqref{Xlin}. Theorem \ref{main} specifies these limits
for scaling exponents $\mbgamma = (\gamma_1,\gamma_2,\gamma_3)$ satisfying \eqref{special}. The general case  $\mbgamma \in \R^3_+$
is treated in Corollary \ref{main}.

\begin{thm} \label{main}  Let $X $ be a linear RF in \eqref{Xlin} with standardized i.i.d. innovations
$\{\vep, \vep(\mbs); \mbs \in \Z^3 \}, \E \vep = 0, \E \vep^2 =1 $ and moving-average coefficients
$a (\mbt)$ in \eqref{aform}, where $\nu >0, q_i>0, c_i >0, i=1,2,3 $  and $\mbq = (q_1.q_2,q_3)$ satisfy
\eqref{qcond}. Moreover, we assume $\lim_{|\mbt| \to \infty} g(\mbt) = 1$ w.l.g.

\smallskip

\noi (i) Let $\frac{1}{2q_{1}} + \sum_{i=2}^3 \frac{1}{q_{i}} < 1$ and $\gamma_{1} q_{1} < \gamma_{2} q_{2}\le \gamma_3 q_3.$
Then
\begin{eqnarray}\label{Fdd1}
\lambda^{-H_1({\small \mbgamma, \mbq})}S^X_{\lambda, {\small \mbgamma}}(\mbx)
&\limfdd&{\cal Y}_1(\mbx).  % B_{{\cal H}_1, 1/2, 1/2}(\mbx)
\end{eqnarray}

\smallskip

\noi (ii) Let $\sum_{i=1}^2 \frac{1}{2q_{i}} + \frac{1}{q_3} < 1 <
\frac{1}{2q_{1}} + \sum_{i=2}^3 \frac{1}{q_{i}}$ and
$\gamma_{1} q_{1} < \gamma_{2} q_{2} < \gamma_{3} q_{3}.$   Then
\begin{eqnarray}\label{Fdd2}
\la^{-H_2({\small \mbgamma, \mbq})}S^X_{\lambda, {\small \mbgamma}}(\mbx)
&\limfdd&{\cal Y}_2(\mbx).  % B_{{\cal H}_1, 1/2, 1/2}(\mbx)
\end{eqnarray}

\smallskip

\noi (iii) Let $1 < \sum_{i=1}^2 \frac{1}{2q_{i}} + \frac{1}{q_{3}}$ and
$\gamma_1 q_1 \le \gamma_{2} q_{2} < \gamma_{3} q_{3}. $  Then
\begin{eqnarray}\label{Fdd3}
\la^{-H_3({\small \mbgamma, \mbq})}S^X_{\lambda, {\small \mbgamma}}(\mbx)
&\limfdd&{\cal Y}_3(\mbx).  % B_{{\cal H}_1, 1/2, 1/2}(\mbx)
\end{eqnarray}

\smallskip

\noi (iv) Let $\sum_{i=1}^2 \frac{1}{2q_{i}} + \frac{1}{q_{3}}  < 1$  and
$\gamma_{1} q_{1} = \gamma_{2} q_{2} < \gamma_{3} q_{3}. $ Then
\begin{eqnarray}\label{Fdd4}
\la^{-H_{12}({\small \mbgamma, \mbq})}S^X_{\la, {\small \mbgamma}}(\mbx)
&\limfdd&{\cal Y}_{12}(\mbx).  % B_{{\cal H}_1, 1/2, 1/2}(\mbx)
\end{eqnarray}

\smallskip

\noi (v) Let  $1 < \frac{1}{2q_{1}} + \sum_{i=2}^3 \frac{1}{2q_{i}}$  and
$\gamma_{1} q_{1} < \gamma_{2} q_{2} = \gamma_{3} q_{3}. $  Then
\begin{eqnarray}\label{Fdd5}
\la^{-H_{23}({\small \mbgamma, \mbq})}S^X_{\la, {\small \mbgamma}}(\mbx)
&\limfdd&{\cal Y}_{23}(\mbx).  % B_{{\cal H}_1, 1/2, 1/2}(\mbx)
\end{eqnarray}

\smallskip

\noi (vi) Let $\gamma_1 q_1 = \gamma_2 q_2 = \gamma_3 q_3.$ Then
\begin{eqnarray}\label{Fdd6}
\la^{-H_0({\small \mbgamma, \mbq})}S^X_{\la, \smbgamma}(\mbx)
&\limfdd&{\cal Y}_0(\mbx).  % B_{{\cal H}_1, 1/2, 1/2}(\mbx)
\end{eqnarray}
The limit RFs and the normalizing exponents in \eqref{Fdd1}-\eqref{Fdd6} are defined
in \eqref{R1def}-\eqref{R6def} and \eqref{H}, respectively.

\end{thm}

To describe the scaling limits in \eqref{Sn} for general $\mbgamma \in \R^3_+$, we need some notation. Let ${\cal P}_3 $ denote the set of all permutations
$\pi = (\pi(1), \pi(2), \pi(3))  $ of $\{1, 2, 3\}$.
Given a RF $ {\cal Y}(\cdot; \mbq, \mbc) = \{{\cal Y}(\mbx; \mbq, \mbc); \mbx \in \R^3_+\}$
depending on vector parameters $\mbc = (c_1,c_2,c_3),
\mbq = (q_1,q_2,q_3) \in \R^3$, and a permutation $\pi = (\pi(1), \pi(2), \pi(3)) \in {\cal P}_3 $,
define a new RF
${\cal Y}^\pi(\cdot; \mbq, \mbc) = \{{\cal Y}^\pi(\mbx; \mbq, \mbc); \mbx \in \R^3_+\}$  by
$$
{\cal Y}^\pi (\mbx; \mbq, \mbc) \ := \ {\cal Y}(\pi \mbx; \pi \mbq, \pi \mbc)
$$
where $\pi \mby  := (y_{\pi(1)}, y_{\pi(2)}, y_{\pi(3)}), \, \mby = (y_1,y_2,y_3)\in \R^3$. The above definition requires some care since
${\cal Y}$ and  ${\cal Y}^\pi$ need not exist simultaneously.
For example, the existence of RFs ${\cal Y}_1(\mbx; \mbq, \mbc)$ in
\eqref{R1def}  and % exists for  $ \frac{1}{2q_1} + \frac{1}{q_2} + \frac{1}{q_3} < 1 $ while
\begin{eqnarray}
{\cal Y}_1^\pi (\mbx; \mbq, \mbc)
&=&\int_{\R^3}  W(\d \mbu)
\int_{\R^3}
\frac{\1(0< u_i < x_i, i= \pi(2), \pi(3), 0< t_{\pi(1)} < x_{\pi(1)}) \d \mbt}{ (c_{\pi(1)}|t_{\pi(1)}-u_{\pi(1)}|^{\frac{q_{\pi(1)}}{\nu}} + \sum_{i=2}^3
c_{\pi(i)}|t_{\pi(i)}|^{\frac{q_{\pi(i)}}{\nu}})^{\nu}}, \nn
\end{eqnarray}
require $ \frac{1}{2q_1} + \frac{1}{q_2} + \frac{1}{q_3} < 1 $ and
$ \frac{1}{2q_{\pi(1)}} + \frac{1}{q_{\pi(2)}} + \frac{1}{q_{\pi(3)}} < 1$, respectively,
and the two conditions are generally different.

From the definition of the partition \eqref{R3G}
it is clear that any $\mbgamma \in \R^3_+$ can be `transformed' into the region \eqref{special} by
a simultaneous permutation of indices of $\gamma_i, q_i$, i.e., for any  $\mbgamma \in \R^3_+$ there exists
a $\pi \in {\cal P}_3$ such that
%that it suffices to describe limit distributions for $\mbgamma \in \R^3_+$ satisfying
\begin{equation}\label{gammaperm}
\gamma_{\pi(1)} q_{\pi(1)} \le \gamma_{\pi(2)} q_{\pi(2)} \le \gamma_{\pi(3)} q_{\pi(3)}
\end{equation}
In general, the above $\pi$ is not unique, e.g., the `well-balanced' points $\mbgamma \in \Gamma_{000} $
satisfy \eqref{gammaperm} for any  $\pi \in {\cal P}_3$.
%for any permutation $\pi \in {\cal P}_3$.
For example, the region
$\gamma_{2} q_{2} \le \gamma_{3} q_{3} \le \gamma_{1} q_{1} =
\Gamma_{\mm 1 1 1} \cup \Gamma_{0 1 1} \cup \Gamma_{000} \cup \Gamma_{\mm 1 0 1} $ corresponds
to \eqref{gammaperm} and $\pi(1) = 2, \pi(2) = 3, \pi(1)=2$.

\begin{cor} \label{main1} Let RF $X$ satisfy the conditions of Theorem \ref{main}.
Let $\pi \in {\cal P}_3$ and
$\mbgamma \in \R^3_+$ satisfy condition \eqref{gammaperm}.

\smallskip

\noi (i) Let $\frac{1}{2q_{\pi(1)}} + \sum_{i=2}^3 \frac{1}{q_{\pi(i)}} < 1$ and $\gamma_{\pi(1)} q_{\pi(1)} < \gamma_{\pi(2)} q_{\pi(2)} 
\le \gamma_{\pi(3)} q_{\pi(3)}.$
Then 
\begin{eqnarray}\label{fdd1}
\la^{-H^\pi_1({\small \mbgamma, \mbq})}S^X_{\la, {\small \mbgamma}}(\mbx)
&\limfdd&{\cal Y}^\pi_1(\mbx; \mbq, \mbc).  % B_{{\cal H}_1, 1/2, 1/2}(\mbx)
\end{eqnarray}

\smallskip

\noi (ii) Let $\sum_{i=1}^2 \frac{1}{2q_{\pi(i)}} + \frac{1}{q_{\pi(3)}} < 1 <
\frac{1}{2q_{\pi(1)}} + \sum_{i=2}^3 \frac{1}{q_{\pi(i)}}$ and
$\gamma_{\pi(1)} q_{\pi(1)} < \gamma_{\pi(2)} q_{\pi(2)} < \gamma_{\pi(3)} q_{\pi(3)}.$   Then
\begin{eqnarray}\label{fdd2}
\la^{-H^\pi_2({\small \mbgamma, \mbq})}S^X_{\la, {\small \mbgamma}}(\mbx)
&\limfdd&{\cal Y}^\pi_2(\mbx; \mbq, \mbc).  % B_{{\cal H}_1, 1/2, 1/2}(\mbx)
\end{eqnarray}

\smallskip

\noi (iii) Let $1 < \sum_{i=1}^2 \frac{1}{2q_{\pi(i)}} + \frac{1}{q_{\pi(3)}}$ and
$\gamma_{\pi(2)} q_{\pi(2)} < \gamma_{\pi(3)} q_{\pi(3)}. $  Then
\begin{eqnarray}\label{fdd3}
\la^{-H^\pi_3({\small \mbgamma, \mbq})}S^X_{\la, {\small \mbgamma}}(\mbx)
&\limfdd&{\cal Y}^\pi_3(\mbx; \mbq, \mbc).  % B_{{\cal H}_1, 1/2, 1/2}(\mbx)
\end{eqnarray}

\smallskip

\noi (iv) Let $\sum_{i=1}^2 \frac{1}{2q_{\pi(i)}} + \frac{1}{q_{\pi(3)}}  < 1$  and
$\gamma_{\pi(1)} q_{\pi(1)} = \gamma_{\pi(2)} q_{\pi(2)} < \gamma_{\pi(3)} q_{\pi(3)}. $ Then
\begin{eqnarray}\label{fdd4}
\la^{-H^\pi_{12}({\small \mbgamma, \mbq})}S^X_{\la, {\small \mbgamma}}(\mbx)
&\limfdd&{\cal Y}^\pi_{12}(\mbx; \mbq, \mbc).  % B_{{\cal H}_1, 1/2, 1/2}(\mbx)
\end{eqnarray}

\smallskip

\noi (v) Let  $1 < \frac{1}{2q_{\pi(1)}} + \sum_{i=2}^3 \frac{1}{2q_{\pi(i)}}$  and
$\gamma_{\pi(1)} q_{\pi(1)} < \gamma_{\pi(2)} q_{\pi(2)} = \gamma_{\pi(3)} q_{\pi(3)}. $  Then
\begin{eqnarray}\label{fdd5}
\la^{-H^\pi_{23}({\small \mbgamma, \mbq})}S^X_{\la, {\small \mbgamma}}(\mbx)
&\limfdd&{\cal Y}^\pi_{23}(\mbx; \mbq, \mbc).  % B_{{\cal H}_1, 1/2, 1/2}(\mbx)
\end{eqnarray}

%\smallskip

%\noi (vi) Let $\gamma_1 q_1 = \gamma_2 q_2 = \gamma_3 q_3.$ Then
%\begin{eqnarray}\label{fdd6}
%p^{-H_0({\small \mbgamma, \mbq})}S_{p {\small \mbgamma}}(\mbx)
%&\limfdd&{\cal Y}_0(\mbx; \mbq, \mba).  % B_{{\cal H}_1, 1/2, 1/2}(\mbx)
%\end{eqnarray}

\end{cor}

The last corollary specifies the scaling limits in the `isotropic' case $q_1 = q_2 = q_3$.

\begin{cor} {\rm Let $X$ satisfy the conditions in Theorem \ref{main}, and
%Consider the `isotropic' case
$q_1 = q_2 = q_3 =: q$.

\smallskip

\noi (I) Let $5/2 < q < 3 $ (Region I). Then
\begin{eqnarray*}
\la^{-H_1({\smbgamma},q)}S^X_{\la, {\small \mbgamma}}(\mbx)
&\limfdd&\sigma_1 B_{{\cal H}_1,1/2,1/2}(\mbx),  \qquad \gamma_1 < \gamma_2 \le \gamma_3,  \\% B_{{\cal H}_1, 1/2, 1/2}(\mbx)
\la^{-H_{12}({\smbgamma},q)}S^X_{\la, {\small \mbgamma}}(\mbx)
&\limfdd&{\cal Y}_{12}(\mbx),  \hskip2.4cm \gamma_1 = \gamma_2 < \gamma_3,  \\
\la^{-H_0({\small \mbgamma},q}S^X_{\la, {\small \mbgamma}}(\mbx)
&\limfdd&{\cal Y}_0(\mbx),  \hskip2.6cm  \gamma_1 = \gamma_2 = \gamma_3.
\end{eqnarray*}

\smallskip

\noi (II) Let $2 < q < 5/2 $ (Region II). Then
\begin{eqnarray*}
\la^{-H_2({\small \mbgamma},q)}S^X_{\la, {\small \mbgamma}}(\mbx)
&\limfdd&\sigma_2 B_{1,{\cal H}_2,1/2}(\mbx),  \qquad \gamma_1 < \gamma_2 < \gamma_3,  \\% B_{{\cal H}_1, 1/2, 1/2}(\mbx)
\la^{-H_{12}({\small \mbgamma},q)}S^X_{\la, {\small \mbgamma}}(\mbx)
&\limfdd&{\cal Y}_{12}(\mbx),  \hskip2.1cm  \gamma_1 = \gamma_2 < \gamma_3,  \\
\la^{-H_{23}({\small \mbgamma},q)}S^X_{\la, {\small \mbgamma}}(\mbx)
&\limfdd&{\cal Y}_{23}(\mbx),  \hskip2.1cm  \gamma_1 < \gamma_2 = \gamma_3,  \\
\la^{-H_0({\small \mbgamma},q)}S^X_{\la, {\small \mbgamma}}(\mbx)
&\limfdd&{\cal Y}_0(\mbx),  \hskip2.3cm  \gamma_1 = \gamma_2 = \gamma_3.
\end{eqnarray*}

\smallskip

\noi (III) Let $3/2 < q < 2 $ (Region III). Then
\begin{eqnarray*}
\la^{-H_3({\small \mbgamma},q)}S^X_{\la, {\small \mbgamma}}(\mbx)
&\limfdd&\sigma_3 B_{1,1,{\cal H}_3}(\mbx),  \qquad \gamma_1 \le \gamma_2 < \gamma_3,  \\% B_{{\cal H}_1, 1/2, 1/2}(\mbx)
\la^{-H_{23}({\small \mbgamma},q)}S^X_{\la, {\small \mbgamma}}(\mbx)
&\limfdd&{\cal Y}_{23}(\mbx),  \hskip1.75cm  \gamma_1 < \gamma_2 = \gamma_3,  \\
\la^{-H_0({\small \mbgamma},q)}S^X_{\la, {\small \mbgamma}}(\mbx)
&\limfdd&{\cal Y}_0(\mbx),  \hskip1.9cm  \gamma_1 = \gamma_2 = \gamma_3.
\end{eqnarray*}
Here, the normalizations and the limit RFs are given as in Theorem \ref{main},
${\cal H}_1 = \frac{7}{2}-q,
{\cal H}_2 = 3-q, {\cal H}_3 =  \frac{5}{2} - q$.
%It follows that the RF $X$ exhibits  complete scaling transition in
%Region II and incomplete scaling transition in Regions I and III.

}
\end{cor}

\begin{remark} We expect that the results of Theorem \ref{main} can be extended to the boundary 
situations 
\begin{eqnarray} \label{H12}
&\frac{1}{2q_1} + \sum_{j=2}^3 \frac{1}{q_j} = 1
\end{eqnarray}
(the boundary between Regions I and II), and
\begin{eqnarray} \label{H23}
&\sum_{j=1}^2 \frac{1}{2q_j} + \frac{1}{q_2} = 1
\end{eqnarray}
(the boundary between Regions II and III), possibly under additional logarithmic normalization. See  also 
\cite{pils2017}, Remark 3.2.     
Note the exponents in \eqref{calH} trivialize in the above cases: 
${\cal H}_1 = 1, {\cal H}_2 = 1/2 $ when \eqref{H12} holds, and ${\cal H}_2 = 1, {\cal H}_3 = 1/2 $
when \eqref{H23} holds. If the above conjecture is true, we can expect  in the limit 
\eqref{Fdd1} and \eqref{Fdd1} a (multiple of) FBS $B_{1,1/2,1/2}$  under \eqref{H12}, and 
a (multiple of) FBS $B_{1,1,1/2}$  under \eqref{H23}, in other words, 
an FBS with {\it all } its Hurst indices equal to 1 and/or 1/2.

\end{remark}

\section{Proof of Theorem \ref{main}}

The proof of Theorem \ref{main} reduces to the central limit theorem for linear forms in i.i.d.r.v.s
$\{\vep(\mbs), \mbs \in \Z^3\} $. Moreover, the limits are written as stochastic integrals w.r.t.
white noise $W (\d\mbu)$ on $\R^3 $. The proof of such limit theorems is facilitated by  the following criterion
generalizing (\cite{book2012},
Prop.14.3.2) to linear forms
\begin{equation} \label{Q}
S(h) := \sum_{\smbs \in \Z^3} h(\mbs) \vep(\mbs)
\end{equation}
with real coefficients  $\sum_{\smbs \in \Z^3} h(\mbs)^2 < \infty$.

\begin{proposition} \label{disc} Let $S(h_\la), \la >0$ be as in \eqref{Q}. Suppose
$h_\la(\mbu)$ are such that for a real-valued function
$f \in L^2(\R^3)$ and some integers $m_i = m_i(\la) \to \infty, \la \to \infty, i=1,2,3 $
the functions
\begin{equation} \label{tildeh}
\tilde h_\la({\mbu}) :=  (m_1 m_2 m_3)^{1/2} h_\la (\lceil m_1u_1\rceil, \lceil m_2u_2\rceil , \lceil m_3 u_3\rceil), \quad
\mbu = (u_1,u_2,u_3) \in \R^3
\end{equation}
tend to $f$ in $L^2(\R^3)$, viz.,
\begin{equation}\label{L2conv}
\|\tilde h_\la - f\|^2 = \int_{\R^3} |\tilde h_\la(\mbu) - f(\mbu)|^2 \d \mbu \to 0, \qquad \la \to \infty.
\end{equation}
Then
\begin{equation}
S(h_\la) \limd I(f) := \int_{\R^3} f(\mbu) W(\d \mbu), \qquad \la \to \infty.
\end{equation}

\end{proposition}

By Cram\'er-Wold device, the proof of finite-dimensional convergence in \eqref{partS} reduces
to the convergence of (scalar) linear combinations $A^{-1}_{\la, \smbgamma} \sum_{k=1}^p \theta_k  S^X_{\la,\smbgamma}(\mbx_k), $ for
any $p \ge 1, \mbx_k \in \R^3_+, \theta_k \in \R, k=1, \cdots, p$ which can be written as linear forms as in \eqref{Q} with
a suitable $h $.
%it suffices to discuss the case $\pi(i) = i, i=1,2,3 $ of the identical permutation $\pi$ only.
For notational convenience, 
we restrict the proof of the last fact to 
%of \eqref{Fdd1}-\eqref{Fdd6} 
to the case $p=1 = \theta_1, \mbx_1 = \mbx $, or to 
the one-dimensional convergence in \eqref{Fdd1}-\eqref{Fdd6} 
since the proof of finite-dimensional convergence is analogous.  Moreover,  for the same reason
we will assume that $\nu = 1, g(\mbt) \equiv 1 $ in \eqref{aform}. We also use the notation 
$V_{\la} (\mbx)$  for normalized sums on the l.h.s. of \eqref{Fdd1}-\eqref{Fdd6}, and drop
$\mbgamma, \mbq $ in the notation of the exponents $H_1(\mbgamma,\mbq), \cdots, H_0(\mbgamma, \mbq)$
in \eqref{H}.

\medskip

\noi {\it Proof of \eqref{Fdd1}.}
Using Proposition \ref{disc}, let
\begin{eqnarray*}
m_i := \lceil \la^{\gamma_i}\rceil,  \quad
\tilde m_i := \la^{\gamma_1 q_1/q_i}, \ i=1,2,3, \quad
\kappa_\la := \frac{(m_1 m_2 m_3)^{1/2} \tilde m_1 \tilde m_2 \tilde m_3}
{\la^{H_1} m_1^{q_1}} \to 1.
\end{eqnarray*}
Then $V_{\la} (\mbx) = S(h_\la)$, where $H_1 = \frac{3\gamma_1+ \gamma_2 + \gamma_3}{2} +
\gamma_1 q_1 (\frac{1}{q_2} + \frac{1}{q_3} - 1)$, see \eqref{H},
\begin{eqnarray}\label{1hla}
h_\la(\mbs)&:=&\la^{-H_1}\sum_{
1\le t_i \le \lfloor \la^{\gamma_i}x_i\rfloor, i=1,2,3} a(\mbt-\mbs) \\
&=&\la^{-H_1}\int_0^{\lfloor \la^{\gamma_1} x_1 \rfloor} \int_0^{\lfloor \la^{\gamma_2} x_2\rfloor}
\int_0^{\lfloor \la^{\gamma_3} x_3\rfloor}
%\int_{\prod_{i=1}^3 (0, \lfloor \la^{\gamma_i} x_i\rfloor]}
%\int_0^{\lfloor \la^{\gamma_3} x_3\rfloor}
\frac{\d t_1 \d t_2 \d t_3}{
c_1|\lceil t_1 \rceil - s_1|_+^{q_1} + \sum_{i=2}^3 c_i|\lceil t_i \rceil - s_i|_+^{q_i }}
\end{eqnarray}
and
\begin{eqnarray*}
&&\tilde h_\la ({\mbu})\ =\  (m_1 m_2 m_3)^{1/2} h_\la (\lceil m_1u_1\rceil, \lceil m_2u_2\rceil, \lceil m_3 u_3\rceil) \\
&&=\ \frac{(m_1 m_2 m_3)^{1/2}}
{\la^{H_1}} \int_0^{\lfloor \la^{\gamma_1} x_1 \rfloor} \int_0^{\lfloor \la^{\gamma_2} x_2\rfloor}
\int_0^{\lfloor \la^{\gamma_3} x_3\rfloor} \frac{\d t_1 \d t_2 \d t_3}{
( c_1|\lceil t_1 \rceil - \lceil m_1 u_1 \rceil |_+^{q_1} + \sum_{i=2}^3 c_i |\lceil t_i \rceil - \lceil m_i u_i \rceil |_+^{q_i} } \\
&&=\ %\frac{(m_1 m_2 m_3)^{1/2} n_1 n_2 n_3} {n^{H(\mbgamma)} n_2^{q_2}}
\kappa_\la
\int_0^{\frac{\lfloor \la^{\gamma_1} x_1 \rfloor}{\la^{\gamma_1}}} \int_0^{\frac{\lfloor \la^{\gamma_2} x_2\rfloor}{\tilde m_2}}
\int_0^{\frac{\lfloor \la^{\gamma_3} x_3\rfloor}{\tilde m_3}}
%&&\
\frac{\d t_1 \d t_2 \d t_3}{
c_1 (\frac{|\lceil \la \tilde m_1  t_1 \rceil - \lceil m_1 u_1 \rceil |_+}{\tilde m_1})^{q_1} + \sum_{i=2}^3 c_i
(\frac{|\lceil \tilde m_i t_i \rceil - \lceil m_i u_i \rceil |_+}{\tilde m_i})^{q_i}
%+(\frac{|\lceil \tilde m_3 t_3 \rceil -  \lceil m_3 y_3 \rceil |_+}{\tilde m_3})^{2q_3/q_1})^{q_1/2}
} \\
&&=\ \int_0^{\lfloor \la^{\gamma_1} x_1\rfloor/\la^{\gamma_1}} \int_\R \int_\R
G_\la(\mbt; \mbu) \d t_1 \d t_2 \d t_3,
\end{eqnarray*}
where
\begin{eqnarray*}
 G_\la(\mbt; \mbu)&:=&
\frac{\kappa_\la
\1\big( -\frac{\lceil m_i u_i\rceil }{\tilde m_i}
< t_i < \frac{ \lceil \la^{\gamma_i}x_i \rceil - \lceil m_i u_i\rceil }{\tilde m_i}, i=2,3 \big)}
{ c_1(\frac{|\lceil \tilde m_1 t_1 \rceil - \lceil m_1 u_1 \rceil |_+}{\tilde m_1})^{q_1} + \sum_{i=2}^3
c_i (\frac{|\lceil \tilde m_i t_i \rceil|_+}{\tilde m_i})^{q_i}
%+ (\frac{|\lceil \tilde m_3 t_3 \rceil\rceil |_+}{\tilde m_3})^{2q_3/q_1})^{q_1/2}
}.
\end{eqnarray*}
Since $\la^{\gamma_i}/m_i \to 1$ and  $m_i/\tilde m_i \to \infty$ due to
$\gamma_i q_i/\gamma_1 q_1 > 1, i=2,3$, we see that
\begin{eqnarray}\label{Gnconv1}
 G_\la (\mbt; \mbu)
&\to&G_1(\mbt; \mbu) \ := \
\frac{\1(0< u_i < x_i, i=2,3)}
{c_1|t_1-u_1|^{q_1} + \sum_{i=2}^3 c_i |t_i|^{q_i}
%+ |t_3|^{2q_3/q_1}  \big)^{q_1/2}
}
\end{eqnarray}
point-wise for any fixed $\mbf u = (u_1,u_2,u_3)\in \R^3,  \mbf t = (t_1,t_2,t_3)\in \R^3, t_1 \ne y_1,
u_i \ne 0, y_i \ne x_i, i=2,3$.
We claim that
\begin{eqnarray}\label{Gla1}
\tilde h_\la({\mbu})
&\to&
\int_{(0, x_1]\times \R^2}
G_1(\mbt; \mbu) \d \mbt  \\
&:=&\1(0< u_i < x_i, i=2,3) \int_{(0, x_1]\times \R^2}  \frac{\d t_1 \d t_2 \d t_3}
{c_1|t_1-u_1|^{q_1} + \sum_{i=2}^3 |t_i|^{q_i}}  \ =: \  f_1({\mbu}) \nn
\end{eqnarray}
point-wise and in $L^2(\R^3)$. Since ${\cal Y}_1(\mbx) = \int_{\R^3}  f_1({\mbu}) W(\d {\mbu})$,
the one-dimensional convergence in \eqref{Fdd1}  follows from \eqref{Gla1}  and Proposition \ref{disc}.

To justify \eqref{Gla1}, note that for all $\lambda \ge 1$ 
\begin{eqnarray}
&&\frac{|\lceil \tilde m_1 t_1 \rceil - \lceil m_1 u_1 \rceil |_+}{\tilde m_1} \ \ge \ \frac{|t_1-u_1|}{2},
\qquad \frac{|\lceil \tilde m_i t_i \rceil |_+}{\tilde m_i}\ \ge \ \frac{|t_i|}{2}, \ i=2,3,\label{dom21}\\
&&\mbox{$\1\big( -\frac{\lceil m_i u_i\rceil }{\tilde m_i}
< t_i < \frac{ \lceil \la^{\gamma_i}x_i \rceil - \lceil m_i u_i\rceil }{\tilde m_i}  \big)
\ \le \ \1(-2 < u_i < x_i + 2) + \1(u_i \le -2, \ \frac{m_i |u_i|}{\tilde m_i} < t_i < \frac{m_i (|u_i| + x_i)}{\tilde m_i})$}
 \nn  \\
&& \hskip6cm + \
\mbox{$\1(u_i \ge x_i+2, \ \frac{m_i (u_i-x_i)}{\tilde m_i} < t_i < \frac{m_i u_i}{\tilde m_i}), \ \ i=2,3.$ }
 \label{dom31}
\end{eqnarray}
Split $\tilde h_\la({\mbu}) = \sum_{j=0}^1 \tilde h_{\la,j}({\mbu})$, where
$\tilde h_{\la, j}({\mbu}) := \int_0^{\lfloor \la^{\gamma_1} x_1\rfloor/\lambda^{\gamma_1}} \int_\R \int_\R
G_{\la, j}(\mbt; \mbu) \d \mbt $ and
\begin{eqnarray*}
G_{\la, 0}(\mbt; \mbu)&:=&G_{\la}(\mbt; \mbu) \1(
-2 < u_i < x_i + 2, i=2,3), \\ % \quad \lambda_1 := (\gamma_2 q_2/q_1) - 1 >0, \\
%G_{n1}(\mbt; \mby)&:=&G_{n}(\mbt; \mby) \1(|y_1| \le 2n^{-\lambda},
%-2 < y_3 < x_3 + 2, ),  \\
G_{\la,1}(\mbt; \mbu)&:=&G_{\la}(\mbt; \mbu) \1(u_i \not\in  (-2,x_i + 2)\  (\exists i =2,3)).
%G_{n3}(\mbt; \mby)&:=&G_{n}(\mbt; \mby) \1(|y_1| \le 2n^{-\lambda_1}, y_3 \not\in  (-2,x_3 + 2)).
\end{eqnarray*}
Relation \eqref{Gla1} follows from
\begin{equation}\label{L2conv11}
\|\tilde h_{\la,0} - f_1\|^2  \to 0  \qquad \text{and} \qquad
\|\tilde h_{\la, 1} \|^2 \to 0.
\end{equation}
The first relation in \eqref{L2conv11} follows from
\eqref{Gnconv} and the dominated convergence theorem since \eqref{dom21}-\eqref{dom31} imply
the dominating bound
\begin{eqnarray*}
0 \le \tilde h_{\la,0}({\mbu})&\le&C \1(-2 < u_i < x_i + 2, i=2,3)  \int_{(0, x_1]\times \R^2}%  \int_0^{x_1} \int_\R \int_\R
\frac{\d \mbt}
{|t_1-u_1|^{q_1} + \sum_{i=2}^3 |t_i|^{q_i}} \ =: \   \bar h({\mbu})
\end{eqnarray*}
with $\bar h \in  L^2(\R^3)$; see the proof of Theorem \ref{exist} (i1).
With $\rho_i
:= \gamma_i - \gamma_1 q_1/q_i >0, i=2,3$,
the second relation in \eqref{L2conv11} follows from
\begin{eqnarray*}
\|\tilde h_{\la, 1} \|^2&\le&C\int_{\R^3} \d \mbu
\Big\{\int_{(0,x_1] \times \R^2} %\int_\R \int_\R
G_1(\mbt; \mbu) \sum_{i=2}^3 (\1 (|t_i| > |u_i| \la^{\rho_i})  \d \mbt \Big\}^2 \ = \  o(1)
\end{eqnarray*}
since $\int_{(0,1] \times \R^2} %\int_\R \int_\R
G_1(\mbt; \cdot )  \d \mbt \in L^2(\R^3)$ and
$\sum_{i=2}^3 \1 (|t_i| > |u_i| \la^{\rho_i}) \to 0 $ for any $\mbt, \mbu \in \R^3$ fixed.
This proves \eqref{L2conv11} and completes the proof of \eqref{Fdd1}.

\medskip

\noi {\it Proof of \eqref{Fdd2}.}
We use Proposition \ref{disc} with
\begin{eqnarray}\label{m2}
&m_1 :=   [\la^{\gamma_2 q_2/q_1}], \quad  m_2:= \lceil \la^{\gamma_2}\rceil, \quad  m_3 := \lceil \la^{\gamma_3}\rceil, \quad
\tilde m_1 := \la^{\gamma_1}, \quad \tilde m_2 := \la^{\gamma_2},  \\
&\tilde m_3 := m_2^{q_2/q_3} \sim \la^{\gamma_2 q_2/q_3}, \qquad
\kappa_\la := \frac{(m_1 m_2 m_3)^{1/2} \tilde m_1 \tilde m_2 \tilde m_3}
{\la^{H_2} m_2^{q_2}} \to 1. \nn
\end{eqnarray}
Then $V_{{\small \mbgamma}} (\mbx) = S(h_\la)$, where $h_\la $ is defined as in \eqref{1hla} with $H_1$ replaced by $H_2$,
and for $\tilde h_\la({\mbu}) $ defined in  \eqref{tildeh}  we have the integral representation
%\begin{eqnarray*}
$\tilde h_\la({\mbu}) = \int_{\prod_{i=1}^2 (0, \lfloor \lambda^{\gamma_1} x_i\rfloor/\tilde m_i]\times \R}  
G_\la(\mbt; \mbu) \d \mbt, $ where
\begin{eqnarray}\label{Gla2} 
 G_\la(\mbt; \mbu)&:=&
\frac{\kappa_\la
\1\big( -\frac{\lceil m_3 u_3\rceil }{\tilde m_3}
< t_3 < \frac{ \lceil \la^{\gamma_3}x_3 \rceil - \lceil m_3 u_3\rceil }{\tilde m_3}  \big)}
{  c_1(\frac{|\lceil \tilde m_1 t_1 \rceil - \lceil m_1 u_1 \rceil |_+}{m_2^{q_2/q_1}})^{q_1} +
c_2 (\frac{|\lceil m_2 t_2 \rceil - \lceil m_2 u_2 \rceil |_+}{m_2})^{q_2} +
c_3 (\frac{|\lceil \tilde m_3 t_3 \rceil\rceil |_+}{m_2^{q_2/q_3}})
^{q_3}}.
\end{eqnarray}
Using $\gamma_2 q_2/\gamma_1 q_1 > 1, \gamma_3 > \gamma_2 q_2/q_3$  we see that
$m_3/\tilde m_3 \to \infty,  m_3/\la^{\gamma_3} \to 1, \la^{\gamma_1} /m_2^{q_2/q_1} \to 0, m_1/m_2^{q_2/q_1} \to 1 $
%m_2/n_2 \to 1,
%\tilde m_3/ m_2^{q_2/q_3} \to 1 $
and hence
\begin{eqnarray}\label{Gnconv}
 G_\la(\mbt; \mbu)
&\to&G_2(\mbt; \mbu) \ := \
\frac{\1(0< u_3 < x_3 )}
{c_1|u_1|^{q_1} + c_2|t_2 -u_2|^{q_2}
+ c_3 |t_3|^{q_3}}
\end{eqnarray}
point-wise for any fixed $\mbf u = (u_1,u_2,u_3)\in \R^3,  \mbf t = (t_1,t_2,t_3)\in \R^3, u_1 \ne 0, u_2 \ne t_2,
t_3 \ne 0, u_3 \ne 0, u_3 \ne x_3$.
%due to $\gamma_2 q_2/q_1 > 1, \gamma_3 > \gamma_2 q_2/q_3$.
Note %Therefore since
$G_2(\mbt; \mby)$ in \eqref{Gnconv}  does not depend on $t_1$. Then
%we can expect that
\begin{eqnarray}\label{Gla2}
\tilde h_\la({\mbu})
&\to&
\int_{(0, x_1]\times (0, x_2]\times \R}
G_2(\mbt; \mbu) \d \mbt \\
&=&x_1 \1(0< u_3 < x_3) \int_{(0, x_2] \times \R}
\big(c_1|u_1|^{q_1} + c_2|t_2 -u_2|^{q_2}
+ c_3|t_3|^{q_3}\big)^{-1} \d t_2 \d t_3
 \ =: \  f_2({\mbu}) \nn
\end{eqnarray}
point-wise and in $L^2(\R^3)$. Since ${\cal Y}_2(\mbx) = \int_{\R^3}  f_2({\mbu}) W(\d {\mbu})$,
the one-dimensional convergence in \eqref{Fdd2}  follows from \eqref{Gla2}  and Proposition \ref{disc}.
The proof of \eqref{GL2} uses  the dominated convergence and a similar argument as in
\eqref{GL2} and we omit the details.

\medskip

\noi {\it Proof of \eqref{Fdd3}.}  Let
\begin{eqnarray}\label{m3}
m_i :=   \lceil \la^{\gamma_3 q_3/q_i}\rceil, \quad \tilde m_i := \la^{\gamma_i},
\quad i=1,2,3, \quad  %
%m_3 := \lceil n^{\gamma_3}\rceil, \quad
\kappa_\la := \frac{(m_1 m_2 m_3)^{1/2} \tilde m_1 \tilde m_2 \tilde m_3}
{\la^{H_3} m_3^{q_3}} \to 1.
\end{eqnarray}
Then $V_{\la} (\mbx) = S(h_\la)$ and $ \tilde h_\la ({\mbu}) = (m_1 m_2 m_3)^{1/2} h_\la (\lceil m_i u_i\rceil, i=1,2,3)
= \int_{ \prod_{i=1}^3 (0, \lfloor \la^{\gamma_i} x_i\rfloor/\tilde m_i]}
G_\la(\mbt; \mbu) \d \mbt, $  where
\begin{eqnarray}
G_\la(\mbt; \mbu)&:=&\frac{\kappa_\la}{
c_1(\frac{|\lceil \tilde m_1 t_1 \rceil - \lceil m_1 u_1 \rceil |_+}{m_3^{q_3/q_1}})^{q_1} +
c_2 (\frac{|\lceil \tilde m_2 t_2 \rceil - \lceil m_2 u_2 \rceil |_+}{m_3^{q_3/q_2}})^{q_2} +
c_3 (\frac{|\lceil \tilde m_3 t_3 \rceil -  \lceil m_3 u_3 \rceil |_+}{m_3})^{q_3}}\nn \\
&\to&\frac{1}{ c_1 |u_1|^{q_1} + c_2 |u_2|^{q_2} +
c_3 |t_3-u_3|^{q_3}} \ =: \ G_3(\mbt; \mbu) \label{Gla3}
\end{eqnarray}
point-wise for any fixed $\mbu = (u_1,u_2,u_3)\in \R^3,  \mbt = (t_1,t_2,t_3)\in \R^3, u_1 \ne 0, u_2 \ne 0,
t_3 \ne u_3 \ne 0$ in view of $\tilde m_1/m_3^{q_3/q_1} \to 0, \,
m_1/m_3^{q_3/q_1} \to 1,\,  \tilde m_2/m_3^{q_3/q_2} \to 0, \,   m_2/m_3^{q_3/q_2} \to 1 $
which follow from the definitions of $m_i, \tilde m_i, i=1,2,3 $ in \eqref{m3} and the inequalities
$\gamma_3 q_3/q_i >\gamma_i, i=1,2$. Note
$G_3(\mbt; \mbu)$ in \eqref{Gla3} does not depend on $t_i, i=1,2$. Also note 
${\cal Y}_3(\mbx) = \int_{\R^3}  f_3({\mbu}) W(\d {\mbu})$, where 
$f_3(\mbu) := \int_{(0, x_1]\times (0, x_2]\times (0,x_3]} G_3 (\mbt; \mbu) \d \mbt$. 
We omit the details of the proof of the convergence $\tilde h_\la \to  f_3 $ in $L^2(\R^3)$, which are similar 
to those in the proof of   \eqref{Fdd1} and \eqref{Fdd2}.

\medskip

\noi {\it Proof of \eqref{Fdd4}.} Let $m_i, \tilde m_i, i=1,2,3, \kappa_\la $ be defined as in \eqref{m2}.  
Note $H_2 = H_{12} = 3(\gamma_1 + \gamma_2)/2 + \gamma_1 q_1/q_3 - \gamma_1 q_1           $ 
for $\gamma_1 q_1 = \gamma_2 q_2$. Then for $G_\la(\mbt; \mbu)$ defined in \eqref{Gla2} we have 
the point-wise convergence
\begin{eqnarray*} %\label{Gnconv}
G_\la(\mbt; \mbu)
&\to&G_{12}(\mbt; \mbu) \ := \
\frac{\1(0< u_3 < x_3 )}
{\sum_{i=1}^2 c_i|t_i-u_i|^{q_i} 
+ c_3 |t_3|^{q_3}}
\end{eqnarray*}
c.f. \eqref{Gnconv}, for any fixed $\mbf u = (u_1,u_2,u_3)\in \R^3,  \mbf t = (t_1,t_2,t_3)\in \R^3, u_i \ne t_i, i=1,2, 
t_3 \ne 0, u_3 \ne 0, x_3$. Moreover, 
${\cal Y}_{12}(\mbx) = \int_{\R^3}  f_{12}({\mbu}) W(\d {\mbu})$, where
$f_{12} (\mbu) := \int_{(0, x_1]\times (0, x_2]\times (0,x_3]} G_{12} (\mbt; \mbu) \d \mbt$. 
The details of the convergence $\tilde h_\la ({\mbu}) 
:= \int_{ \prod_{i=1}^3 (0, \lfloor \la^{\gamma_i} x_i\rfloor/\tilde m_i]}
G_\la(\mbt; \mbu) \d \mbt \to  f_{12} (\mbu)$ in $L^2(\R^3)$ 
are similar as above.

\medskip

\noi {\it Proof of \eqref{Fdd5}.} Let $m_i, \tilde m_i, i=1,2,3, \kappa_\la $ be defined as in \eqref{m3}.  
Note $H_3 = H_{23} = \gamma_1 + 3(\gamma_2 + \gamma_3)/2 + \gamma_2 q_2/2q_1 - \gamma_2 q_2            $
when $\gamma_2 q_3 = \gamma_3 q_3$. Then for $G_\la(\mbt; \mbu)$ defined in \eqref{Gla3} we have
the point-wise convergence
\begin{eqnarray*}
G_\la(\mbt; \mbu)&\to&\frac{1}{ c_1 |u_1|^{q_1} + 
\sum_{i=2}^3 c_i |t_i-u_i|^{q_i}} \ =: \ G_{23}(\mbt; \mbu)
\end{eqnarray*}
point-wise for any fixed $\mbu = (u_1,u_2,u_3)\in \R^3,  \mbt = (t_1,t_2,t_3)\in \R^3, u_1 \ne 0,   
u_i \ne t_i, i=2,3$.  Moreover,
${\cal Y}_{23}(\mbx) = \int_{\R^3}  f_{23}({\mbu}) W(\d {\mbu})$, where
$f_{23} (\mbu) := \int_{(0, x_1]\times (0, x_2]\times (0,x_3]} G_{23} (\mbt; \mbu) \d \mbt$.
The details of the convergence $\tilde h_\la ({\mbu})
:= \int_{ \prod_{i=1}^3 (0, \lfloor \la^{\gamma_i} x_i\rfloor/\tilde m_i]}
G_\la(\mbt; \mbu) \d \mbt \to  f_{23} (\mbu)$ in $L^2(\R^3)$
are similar and omitted.

\medskip

\noi {\it Proof of \eqref{Fdd6}.}  Let $m_i :=   \lceil \la^{\gamma_i}\rceil,  \tilde m_i := \la^{\gamma_i},
i=1,2,3, \, \kappa_\la := \prod_{i=1}^3 m_i^{1/2} \tilde m_i/\la^{H_0} m_1^{q_1} \to 1.$  
As noted above, in this case $H_0 = \gamma_1 q_1((3/2)\sum_{i=1}^3 1/q_i - 1)$ agrees with any of 
$H_1, \cdots, H_{23}$ in the above proof. We also see that $G_\la(\mbt; \mbu)$ defined in \eqref{Gla1}-\eqref{Gla3} 
tends to  $G_{0}(\mbt; \mbu)$, viz., 
\begin{eqnarray*}
G_\la(\mbt; \mbu)&\to&\frac{1}{ 
\sum_{i=1}^3 c_i |t_i-u_i|^{q_i}} \ =: \ G_{0}(\mbt; \mbu)
\end{eqnarray*}
point-wise for any fixed $\mbu = (u_1,u_2,u_3)\in \R^3,  \mbt = (t_1,t_2,t_3)\in \R^3, 
u_i \ne t_i, i=1,2,3$, and 
${\cal Y}_{0}(\mbx) = \int_{\R^3}  f_{0}({\mbu}) W(\d {\mbu})$, where
$f_{0} (\mbu) := \int_{(0, x_1]\times (0, x_2]\times (0,x_3]} G_{0} (\mbt; \mbu) \d \mbt$.  
We omit the rest of the proof  since it is similar as in the previous cases.  Theorem \ref{main} is proved. 
\hfill $\Box$

%\end{description}

\end{document}